\documentstyle[12pt,a4,epsfig]{article}
\newtheorem{theor}{Theorem}

\begin{document} 
\title{Maps of $p$-gons with a ring of $q$-gons}
\author{Michel Deza\\
CNRS and Ecole Normale Sup\'erieure, Paris \\ \and 
Viatcheslav Grishukhin\\
CEMI, Russian Academy of Sciences, Moscow}
\date{} 

\def\d#1{\epsfig{file=#1.ps,width=4cm}}
\def\da#1{\epsfig{file=#1.ps,width=4.4cm}}
\def\dd#1{\epsfig{file=#1.ps,width=3.6cm}}
\def\ddd#1{\epsfig{file=#1.ps,width=3.3cm}}

\maketitle 

\begin{abstract}
We study 3-valent maps $M_n(p,q)$ consisting of a ring of $n$ $q$-gons whose 
the inner and outer domains are filled by $p$-gons, for $p,q \ge 3$. We 
describe a domain in the space of parameters $p$, $q$, and $n$, for which 
such a map may exist.
With four infinite sequences of maps - prisms $M_p(p \ge 3,4)$,
$M_4(4,q \ge 4)$, $M_4(5,5t+2 \ge 7)$, $M_4(5,5t+3 \ge 8)$, we give 26 
sporadic ones. The maps whose $p$-gons form two paths are first two infinite
sequences and 5 maps: $M_{28}(7,5)$, $M_{12}(6,5)$, $M_{10}(5,6)$, 
$M_{20}(5,7)$, $M_{2}(3,6)$.  

\end{abstract}

Let $M_n(p,q)$ be a planar graph with all vertices of degree 3 having 
only $p$-gonal and $q$-gonal faces such that the $q$-gonal faces form 
a ring $R_n$ of $n$ $q$-gons. We want to know for which $n$, $p$ and $q$ such 
maps exist. Clearly, the dual of $M_{n}(p,q)$ is
a triangulation, in which all vertices are $p$-valent, except
$n$ $q$-valent vertices forming an induced simple circuit. 
 
We dismiss trivial maps $M_{3}(2,2+4t)$ (with $t$ consecutive digons on each
of 3 edges of dual triangle) and $M_{3}(p,p)$, $3 \le p \le 5$, of simple
Platonic polyhedra. So, one can consider only the case $q \ge 4$ since the map 
$M_{3}(3,3)$ is unique $M_{n}(p,3)$.
We represent the ring $R_n$ as a ring of quadrangles as follows. Each 
quadrangle has two pairs of opposite edges. Edges of one pair belong to 
neighboring quadrangles. Edges of another pair belong to the outer and inner 
boundaries of the ring. In order to transform a quadrangle into a $q$-gon, we set 
$q-4$ vertices on the edges of the quadrangle belonging to the outer and inner 
boundaries. If $(i,j)$ is the pair of numbers of the vertices on the outer and 
inner edges of quadrangles, then the following  pairs are possible: 
$(0,q-4)$, $(1,q-5)$,...,$(q-4,0)$. So any of our maps will be a decoration of 
a $n$-gonal prism.

Denote by $I_n$ and $O_n$ the parts of $M_n(p,q)$ consisting of the
inner and outer domains of $R_n$, with the common boundary of $R_n$.

If $q=4$, then the map $M_n(p,4)$ exists only for $n=p$. In this case $I_n$ 
and $O_n$ each is a $p$-gon, and $M_p(p,4)$ is the map of a $p$-gonal 
prism with quadrangle lateral faces. Hence, from now on we suppose that $q \ge 5$. 

Consider the part $I_n$ of $M_n(p,q)$ consisting of $p$-gons in the inner 
domain of $R_n$. Vertices of $I_n$ have degrees 2 and 3. Let $v_2$ and $v_3$
be the numbers of vertices of degree 2 and 3 lying on the common boundary of $I_n$ 
and $R_n$. The vertices of degree 2 are end-vertices of edges common to two 
adjacent $q$-gons. Obviously, there are $v_2=n$ such vertices. Let $x$ be 
the number of interior vertices of $I_n$. 

We are interested especially in maps $M_n(p,q)$ such that $I_n$ and $O_n$ 
each is a path of $p$-gons. We say in this case that $M_n(p,q)$ {\em has 
two paths of $p$-gons}. 

Now we apply results of \cite{DFG} to the map $I_n$. The boundary of $I_n$ 
generates a {\it $p$-gonal boundary sequence} studied in \cite{DFG}. Each 
$p$-gonal sequence has the form $a=a_1a_2...a_k$, where $k$ is the length 
of the sequence $a$. Here $a_i$ is the number of vertices of degree 2 between 
the {\it tails} $t_i$ and $t_{i+1}$. Namely, the tail $t_i$ is an edge of
$I_n$ having exactly one end-vertex (of degree 3) on its boundary.Hence we have 
\begin{equation}
\label{kv}
k=v_3 \mbox{ and }\sum_{i=1}^k a_i=v_2. 
\end{equation}

Note that if $p=n$, then there is a degenerate case, when one of the
domains $I_n$ and $O_n$ or both ($p$-prism) may consists of one $p$-gon. In this
case $v_3=0$, what implies $k=0$. The corresponding $p$-gonal
sequence $a$ is degenerate also.

Similarly, let $v'_2$, $v'_3$, $x'$ be the numbers of the corresponding 
vertices of the outer domain $O_n$. The boundary of $O_n$ generates a 
boundary sequence of length $k'=v'_3$. Of course, we have 
\begin{equation}
\label{IO}
v'_2=v_2=n, \hspace{3mm}v'_3+v_3=(q-4)n. 
\end{equation}
The equalities (\ref{kv}) and (\ref{IO}) imply 
\begin{equation}
\label{kn}
\sum_i^ka_i=n, \hspace{3mm}\sum_i^{k'}a'_i=n, \hspace{3mm}k+k'=(q-4)n.  
\end{equation} 
We call the sequence $a'$ the {\em $q$-complement} of $a$, and vice versa.
A sequence $a$ is called {\em self $q$-complemented} if $a'$
is obtained from $a$ by shifting or/and reversing.

\begin{theor}
For $p,q \ge 4$, a map $M_n(p,q)$ may exist only for the following values 
of the parameters $n$, $p$ and $q$. 

{\bf A} If $q=4$, then $n=p$ and $M_p(p,4)$ is the map of a $p$-prism. 

{\bf B} If $q \ge 5$, then $p \le 7$, and if $p \ge 6$, then $q=5$. Besides

1) if $p=7$, then $n \ge 28$ and unique $M_{28}(7,5)$ has two paths of 
$7$-gons; 

2) if $p=6$, then all existing maps are four $M_{12}(6,5)$ and one of them has 
two paths of hexagons (see Figure 1);  

3) if $p=5$, then  

\hspace{5mm}a) if $q=5$, then there exist only two maps $M_6(5,5)$ and one 
$M_5(5,5)$, all three realizing the dodecahedron; 

\hspace{5mm}b) if $q=6$, then $n \le 10$, and all existing maps are: 
$M_5(5,6)$, two maps $M_6(5,6)$, $M_8(5,6)$ and 
(with two paths of pentagons) $M_{10}(5,6)$ (see Figure 1); 

\hspace{5mm}c) if $q=7$, then $n \le 20$ and all existing maps are (except
of undecided case $17 \le n \le19$) following nine:
$M_4(5,7)$, $M_{10}(5,7)$ and (with two paths of pentagons) $M_{20}(5,7)$
(see Figure 2), and (from \cite{Ha}, see Figure 3) 4 maps $M_{12}(5,7)$,
2 maps $M_{16}(5,7)$;

\hspace{5mm} if $q \ge 8$, then if a map $M_n(5,q)$ exists, it has not paths of 
pentagons, and the following maps exist: $M_2(5,10)$, $M_3(5,8)$ and 
$M_4(5,q)$, for any $q \equiv 2,3 (mod \, 5)$, but not for $q=9$; 

4) if $p=4$, then $n \le 4$ and all existing maps are: $M_2(4,8)$, 
$M_3(4,6)$ and (with two paths 
of $q-3$ quadrangles) $M_4(4,q)$ for any $q \ge 4$; 

5) if $p=3$, then there exist exactly two maps $M_2(3,6)$ and $M_3(3,4)$,
where the last map is the special case $p=3$ of {\bf A}.

\end{theor}
{\bf Proof}. The case {\bf A} was described above. Consider 
the case {\bf B}. 

We use the Euler relation $v-e+f=1$ for numbers of vertices $v$, 
edges $e$ and faces $f$ of $I_n$. We have $v=v_2+v_3+x$, $2e=2v_2+3(v_3+x)$, 
$pf=v_2+2v_3+3x=v_2-v_3+3(v_3+x)$. Substituting these values in the Euler 
relation, one get 

\[2v_2-(p-4)v_3+(6-p)x=2p. \]
A similar equality is valid for $O_n$. Summing them and 
using (\ref{IO}), we get 
\begin{equation}
\label{pqn}
((4-p)(q-4)+4)n+(6-p)(x+x')=4p. 
\end{equation}

Consider at first the case $p \ge 7$. Rewrite (\ref{pqn}) as follows:
\[(p-6)(x+x')=(4-(p-4)(q-4))n-4p. \]
Since $p-6>0$ and $x+x' \ge 0$, one has $4-(p-4)(q-4)>0$. 
For $p \ge 7$, this inequality holds only if $q=5$ and $p=7$. For these values 
$p$ and $q$, the equation (\ref{pqn}) takes the form 
\[x+x'=n-28. \]
So, no map $M_{n}(7,5)$ exists for $n<28$. For $n=28$, $x+x'=0$. Hence if
a map $M_{28}(7,5)$ exists, then $I_{28}$ and $O_{28}$ consist each of
a path of heptagons. In fact, such  map exists and it is unique for
$n=28$.

Now we consider the case $p=6$. In this case (\ref{pqn}) takes the form 
\[(12-2q)n=24. \]
Hence $12-2q>0$, i.e. $q<6$. Since $q \ge 5$, we have $q=5$ and $n=12$. 
There exist four maps $M_{12}(6,5)$, and only one of them has two paths of 
hexagons. (See Proposition 9.1 of \cite{DG}.) 

The case $p=5$ is the most rich case. Now (\ref{pqn}) takes the form 
\begin{equation}
\label{xn}
x+x'=20-(8-q)n.  
\end{equation}

For $q<8$, this equality restricts $n$, since $x+x' \ge 0$.

Now we will use Table 1, which is Table 3 of \cite{DFG} extended by the sequences 
of length $k=9$. For given $n$ and $q$, we inspect pentagonal sequences of 
length $k \le \frac{1}{2}(q-4)n$ and such that $\sum_1^k a_i=n$. Call such 
a sequence {\em feasible}. We seek feasible sequences having pentagonal 
$q$-complements. 
At right of a feasible sequence we give in parentheses the map related 
to this sequence.

Consider 
at first the case $q=5$. The equality (\ref{xn}) gives $n \le 6$. But 
now $M_n(5,5)$ is a partition of the plane into pentagons. There exists 
only one such partition and it is the dodecahedron. 
It can be realized as one map $M_5(5,5)$
with a degenerated sequence $a$, and two maps $M_6(5,5)$, one with
the self 5-complemented sequence 222 and another with the sequence
2121 and its 5-complement 33.

\vspace{3mm}
\begin{center}
\footnotesize
\vbox{{\bf Table 1. Pentagonal sequences of length $k$ at most 9.}
\begin{tabular}{|l|l|llllll|} \hline
$k$ & $n$& \multicolumn{6}{c|}{$n= \sum_1^k a_i$, pentagonal sequences}\\ \hline
2 & 6 & 33 & ($M_6(5,5)$)  & & &&\\ \hline
3 & 6 & 222  &($M_6(5,5)$)&& & & \\ \hline
4 & 6 & 2121 &($M_6(5,5)$) & &&& \\  
  & 7 & 3130 & & & &&\\ \hline 
% & 0 & 00000 & ($M_5(5,5)$)& & &&\\ 
5  & 5 & 11111 &($M_5(5,6)$)&& & & \\ 
  & 6 & 21120 & & & &&\\
  & 7 & 30220 & & & &&\\ \hline 
6 & 2 & 100100 &($M_2(5,10)$)&  & & &  \\ 
  & 3 & 101010 &($M_3(5,8)$)&  & & & \\
  & 4 & 110110 &($M_4(5,7)$)&& & &  \\
  & 5 & 201110 & & & &&\\
  & 6 & 202020 &($M_6(5,6)$)  & 210210 &($M_6(5,6)$)& & \\
  & 7 & 301210 & & & &&\\
  & 8 & 311300 && 310310& & & \\ \hline
7 & 4 & 2001100 & & & && \\ 
  & 5 & 2010200 & & & &&\\
  & 6 & 2101200 & & & && \\
  & 7 & 3011200 & & & &&\\
  & 8 & 3102200 & &3020300& & &  \\ \hline
8 & 4 & 20002000  &($M_4(5,8)$)&  & & &  \\ 
  & 6 & 30011100 & &21002100& & &  \\
  & 7 & 30102100 & &30020200 && &  \\
  & 8 & 31012100 & &30103010 &($M_8(5,6)$)  & 22002200 &  \\
  & 9 & 31113000 & &31103100 & &31013010 & \\ \hline
9 & 5 & 300011000 & &&& & \\
  & 6 & 300102000 & &210012000& & & \\
  & 7 & 301012000 & & & &&\\
  & 8 & 310112000 & &310021100 & &301103000 & 220012100 \\
  & 9 & 311022000 & &310203000 && 310022010 & 300300300 \\
  &10 & 311002201 & & &&& \\ \hline 
\end{tabular}}
\end{center}

\vspace{3mm}
Let $q=6$. The equality (\ref{xn}) takes the form $x+x'=20-2n$. Hence 
$n \le 10$. One has to consider feasible pentagonal sequences of length 
$k \le n$. For $2 \le n \le 10$, we find the following feasible pentagonal 
sequences: 

$n=5$, $k=5$: 11111. It is a self 6-complemented sequence giving $M_5(5,6)$. 

$n=6$: for $k=2,3,4,5$ we find the feasible sequences 33, 222, 2121, 21120, 
respectively. (Note that the 5-complement of 222
is pentagonal giving a map $M_6(5,5)$ of the dodecahedron.) For $k=6$, we find 
two feasible sequences 202020 and 210210 (both are self-complemented), 
giving two distinct maps $M_6(5,6)$. 

For $n=7$ all feasible sequences have no pentagonal 6-complement.  

For $n=8$ there is unique feasible sequence 30103010 of length 8, which 
is self 6-complemented and gives the map $M_8(5,6)$. 

For $n=9$: the 6-complement of unique feasible sequence is not pentagonal.

If $n=10$, then $x+x'=0$. There exists a map $M_{10}(5,6)$ having two
paths of pentagons. It is unique for $n=10$. 

Let $q=7$. The equality (\ref{xn}) takes the form $x+x'=20-n$. Hence 
$n \le 20$. We have to consider feasible sequences of length 
$k \le \frac{3}{2}n$. For $n=2$ and $n=3$ there are no feasible 
sequences. Hence maps $M_n(5,7)$ do not exist  for $n \le 3$ . 

For $n=4$, Table 1  shows that amongst of all sequences of length $k \le 6$ 
there is unique feasible sequence, namely 110110, This sequence is 
self 7-complemented. Hence in this case both maps $I_4$ and $O_4$ are 
isomorphic and contain 8 pentagons each. We obtain the map $M_4(5,7)$, 
starting the sequence $M_4(5,5t+2)$.

For $n=5$, Table 1 shows that amongst all sequences of length $k \le 7$ there 
are 3 feasible sequences 11111, 201110 and 2010200. But 
the 7-complements of all these sequences are not pentagonal. So, there is 
no map $M_5(5,7)$. 

For $n=6$, Table 1 gives 11 feasible sequences, but the 7-complements of all these 
sequences are not pentagonal. 

There exist $M_{10}(5,7)$ with pentagonal sequence $bb$ for $b=3010010$; its
$7$-complement is pentagonal sequence $bb$ for $b=21001001$.

Harmuth ( \cite{Ha}) checked by computer all maps $M_{n}(5,7)$ with $n \le 16$.
He got six new maps (four with $n=12$ and two with $n=16$); see them on 
Figure 3. Amongst them, only 3rd and 4th have $I_n \neq O_n$. The 4th map
$M_{12}(5,7)$ has non-periodic pentagonal sequence, other five have it of form
$bb$.

Note the extremal case $n=20$. In this case $k+k'=60$. There exists a 
feasible self 7-complemented sequence of length 30. This sequence has the 
form $bb$, where $b=300100101011011$. The inner and outer maps $I_{20}$ 
and $O_{20}$ of the map $M_{20}(5,7)$ are isomorphic. 
In this case $x=x'=0$ and  $v_3=v'_3=\frac{3}{2}n=30$. Hence each of them  
consists of a path of $f=1+\frac{1}{2}v_3=1+\frac{1}{2}v'_3=16$ pentagons. 
This map is unique for $n=20$; its symmetry is $S_4$.

If $q \ge 8$, the equality (\ref{xn}) gives no restriction on $n$. 
There exist maps $M_3(5,8)$, $M_2(5,10)$, corresponding
to pentagonal self q-complemented sequences with parameters $(k,n)=(6,3),
(6,2)$, respectively.
Moreover, there exist maps $M_4(5,q)$ for any $q=5t+3 \ge 8$ and
$q=5t+2 \ge 7$. They have self $q$-complementary sequences $bb$ with 
$b=20...0$ ($5t-2$ zeros) and $b=110...0$ ($5t-4$ zeros), respectively.
The first is a decorated (not on the $4$-ring) map $M_4(4,5t+3)$; the 
second is a decorated $M_4(4,5)$. See the smallest case $t=1$ (i.e. $q=7,8$)
of those maps on Figure 2; the larger ones come by repetition of their
$I_n, O_n$ $t$ times. One can check non-existence of $M_4(4,9)$, using that its
sequence $a$ should be of form $bb$ (since it period divides 4) and has $k=10$
(since Table 1 has no needed sequence with $n=4$ and $k \le 9$ or $k' \le 9$).

Let $p=4$. In this case (\ref{pqn}) takes the form 
\[n=4-\frac{1}{2}(x+x'), \]
that implies that $n \le 4$. An inspection shows that there exist the
following maps 
\[M_2(4,8), \hspace{3mm}M_3(4,6), \hspace{3mm}M_4(4,q), \mbox{  }q \ge 4. \]
Since $x+x'=0$ for $n=4$, the inner and outer domains of the map $M_4(4,q)$, 
each consists of a path of $q-3$ quadrangles. 

The map $M_4(4,q)$ can be obtained from the prism $M_p(p,4)$
for $p=q$ as follows. Select a 4-gon $Q$ in the ring $R_q$ of
$M_q(q,4)$ and put in $Q$ new $q-4$ edges parallel to edges of $Q$
common with the two $q$-gons. Then $Q$ will be replaced by a path
of $q-3$ quadrangles. The two 4-gons, which were neighboring to $Q$,
became $q$-gons that form a ring of four $q$-gons together with the
old two $q$-gons. (The dual of $M_4(4,5)$ is the map of {\it snub dishpenoid};
its faces are 12 regular triangles.)

Now, consider the last case $p=3$. In this case (\ref{pqn}) takes the form
\[x+x'=4-\frac{1}{3}qn. \]
It implies that $qn \le 12$ and $qn$ is a multiple of 3. Since $q\ge 4$,
only two pairs $(n,p)=(2,6)$ and $(n,p)=(3,4)$ are possible. Both
corresponding maps $M_2(3,6)$ and $M_3(3,4)$ exist with isomorphic
$I_n$ and $O_n$.
\hfill $\Box$
 
\vspace{1mm}
Clearly, $M_n(p,q)$ has two paths of $p$-gons if and only if $x+x'=0$ and such
paths consist of $ \frac{n-4}{p-4}$ $p$-gons. So, in this case 
(\ref{pqn}) implies $n= \frac{4p}{4-(p-4)(q-4)}$ and $(p,q)=(5,6),(5,7),(6,5),
(7,5),(3,6)$. Any such map comes from $M_4(4, \frac{n-4}{p-4}+3)$ by addition of $n-4$ edges on the $4$-ring.

\vspace{2mm}
{\bf Corollary} {\em Besides two infinite sequences of maps: 
$M_p(p,3)$ and $M_4(4,q)$, there are only 
5 maps $M_n(p,q)$, having two paths of $p$-gons:
$M_{2}(3,6)$, one of 4 maps $M_{12}(6,5)$, $M_{10}(5,6)$,
$M_{28}(7,5)$, $M_{20}(5,7)$ (see Figures 1,2).}$\Box$

{\bf Remarks:}

(i) Maps $M_n(5,6)$ and $M_n(6,5)$ are instances of {\em fullerenes}, 
known in Organic Chemistry; they are 9 maps (of cases {\bf B}2) and
{\bf B}3)b) of Theorem 1) given on Figure 1. In notation $F_m(Aut)$, stressing their number of vertices (carbon atoms) and the maximal symmetry, those 9 maps
are: $F_{36}(D_{2d})$, $F_{44}(D_{3d})$, $F_{48}(D_{6d})$, $F_{44}(D_2)$
and $F_{30}(D_{5h})$, $F_{32}(D_2)$, $F_{40}(D_2)$,
$F_{32}(D_{3d})$, $F_{36}(D_{2d})$.

The 4 simple polyhedra with only 4- and 6-gonal faces and at most 16
vertices, correspond to 4 maps $M_4(4,4), M_6(6,4)$ (4- and 6-prisms) and
$M_3(4,6), M_4(4,6)$. 

The 26 maps $M_n(p,q)$, given here (which are not covered by families 
$M_p(p \ge 3,4)$, $M_4(4,q \ge 4)$, $M_4(5,q \equiv 2,3 (mod \, 5)$) are: 
3 maps of the dodecahedron, 5 small maps $M_n(p,q)$ (with 
$(n; p,q)$=$(2; 3,6)$, $(3; 4,6)$,
$(2; 4,8)$, $(3; 5,8)$,$(2; 5,10)$; see Figure 2), 9 fullerenes and 9 
{\it azulenoids} (a
chemical term for the case $ \{ p,q \} = \{ 5,7 \}$). All undecided maps
$M_n(p,q)$ are either other azulenoids ($M_n(7,5)$,
$n \ge 29$, or $M_n(5,7)$, $17 \le n \le 19$), or, prior to \cite{SM}, other {\it fulleroids}
(a chemical term for the case $p=5$) $M_n(5,q)$ (with
$q \ge 8$, $n \ge 4$).

\cite{SM} gave new constructions of fulleroids $M_n(5,q)$ for 
$n=8, q \equiv 1,4 (mod \, 5), q \ge 9$ (by S.Madaras) and four cases
(by R.Sotak): 1) $n=6, q \equiv 0,1,4 (mod \, 5), q \ge 10$,
2) $n=10, q \equiv 0 (mod \, 10), q \ge 20$, 
3) $n=14+4k, q=10, k \ge 0$,
3) $n=12+4k, q \equiv 2,3 (mod \, 5), q \ge 17, k \ge 0$. So any even $n \ge 0$ and any $q \ge 4$ are realized (separately) by some map $M_n(5,q)$, but 
amongst impossible pairs $(n,q)$ there are, for example, $(4,9)$, $(n,6)$ for
$n>10$, $(n,7)$ for $n>20$ and, using \cite{Ma}, $(n,q \equiv 0 (mod \, 10))$
for odd $n$. All known $M_n(p,q)$ with {\it odd} $n$ are, besides of prisms
$M_n(n,4)$ with odd $n$, the map $M_3(4,6)$ and 3 maps $M_5(5,q)$ with
$q=4,5,6$. There is a $M_n(5,10)$ for each $n \equiv 2 (mod \, 4)$ (except of
undecided case $n=10$) and an
infinity of maps $M_n(5,q)$ for any $n \equiv 0 (mod \, 4)$ and for 
$n=6,10$. The number of $M_n(5,q)$ with fixed $q$ is finite for $q \le 6$, but
it is infinite for $q=10$ and $q \equiv 2,3 (mod \, 5), q \ge 17$.

(ii) Examples of two maps having the
same graph $G$ of $p$-gons in $I_n$, are: $M_3(4,6)$, one of 2 $M_6(5,5)$ (for
$G=K_3$); $M_2(4,8)$, the second $M_6(5,5)$ (for $G=K_4-e$); $M_5(5,6)$,
$M_5(5,5)$ (for $G$ being $5$-wheel); one of two $M_6(5,6)$, one of four $M_{12}(6,5)$ 
(for $G=K_{3 \times 2}-K_3$).

The maps, given here (except of $M_5(5,5)$, one of two $M_6(5,5)$, 
$M_{10}(5,7)$ and 3rd, 4th $M_{12}(5,7)$ from Figure 3) are
{\it self-complementary} (i.e. they have isomorphic domains $I_n$ and $O_n$)
since the corresponding $p$-gonal sequences are self $q$-complemented.
Suppose, without loss of generality, that $I_n$ contain no more $p$-gons than
$O_n$. Then exceptional maps $M_5(5,5)$, $M_6(5,5)$, 
$M_{10}(5,7)$, $M_{12}(5,7)$   have (the
graph of $p$-gons in $I_n$) $G= K_1, K_2$, 
$P_{0,1,...,9}$+ $(1,3)$+$(2,4)$+$(5,7)$+$(6,8)$,
$P_{0,1,...,9}$+ $(0,2)$+$(7,9)$,
$P_{0,1,...,9}$+ $(1,3)$+$(2,4)$+$(7,9)$.

All $p$-gonal sequences of the maps in the paper (except of non-periodical one
for 4th map on Figure 3) are:
either 222222 (for a $M_{12}(6,5)$), or 11111 (for $M_5(5,5)$ and $M_5(5,6)$), or
of form $bbb$ (with $b=$ 1, 10, 2, 20, 31 for 5 maps $M_3(4,6)$, a 
$M_6(5,5)$, $M_3(5,8)$, $M_6(5,6)$, a $M_{12}(6,5)$), or of form $bb$ (including $b=$ 1, 10, 100, 110...0, 20...0 (with any number of zeros), 21, 210, 321, 411, 3010, 30101, 3010010). Clearly, the period divides $n$ and indicates the maximal 
symmetry of corresponding polyhedron.

(iii) Denote by $M_{n}^{k}(p,q)$ {\it $k$-valent} analogs of maps
$M_n(p,q)$, i.e. $k$-valent $(p,q)$-map, such that $q$-gons form a simple circuit.
Unique map $M_{n}^{4}(p,3)$ is the map $M_{2p}^{4}(p,3)$ of $p$-gonal 
antiprism; all $M_{n}^{4}(3,4)$ are 3 maps with $n=4,6,8$, which are 
decorations of
cube, octahedron and rhombic dodecahedron, respectively.

The only $k$-valent planar graphs, admitting
more than one map $M_{n}^{k}(p,p)$, are the dodecahedron and
the icosahedron. The dodecahedron has three maps: with $(n,G)=(5,K_1), (6,K_2)$
and (self-complementary one) with $(6,K_3)$. The icosahedron has six 
maps $M_{n}^{5}(3,3)$:
with $(n,G)=(9,P_1)$, $(10,P_2)$, $(11,P_3)$ and 
(self-complementary ones) with $(12,P_4)$, $(12,K_{1,3})$, $(10,C_5)$.

(iv) Denote by $M_{n_1,...,n_t}(p,q)$ any 3-valent planar map, consisting of
$p$- and $q$-gons only, where all $q$-gons are organized into $t, t \ge 2,$
rings of length $n_1,...,n_t \ge 3$. We found 3 infinite families
(of $M_{3,3}(6,4)$, $M_{6,6}(6,5)$, $M_{3,3,3,3}(6,5)$) and 30 sporadic maps:
3 fullerenes ($M_{3,3}(5,6)$, $M_{5,5}(5,6)$, $M_{3,6}(5,6)$=
$M_{3,9}(6,5)$),
15 azulenoids (12 $M_{n_1,n_2}(5,7)$, $M_{3,3,3,3}(5,7)$,
$M_{3,8,3}(5,7)$, $M_{6,9,3}(5,7)$=$M_{3,15,12}(7,5)$) and 12 maps
$M_{3,...,3}(p,4)$ with (t,p)=

(4,7),(4,8),(4,8),(4,9),(6,10),(8,9),(8,9),(8,12),(12,11),(12,11),(16,13),(20,15). 

(v )For $n=2,3$ only the maps
$M_2(p,2p)$ and $M_3(p,2p-2)$ for $p=3,4,5$ exist. The maps with $n=2$ are not
polyhedral: two $2p$-gons have two edges in common.

The maps $M_2(p,2p)$ with $p=3,4,5$ have symmetry $D_{2h}$; two maps from Figure 1 have
symmetry $D_{3d}$ (namely, fullerenes $M_{12}(6,5)$ and $M_6(5,6)$, which are the first and the last
in the right column). Those five maps are only centrally-symmetric ones amongst the maps of
Figures 1 and 2. They can be folded on real projective plane: the ring of $n$ $q$-gons became
M\"{o}bius band from $\frac{n}{2}$ $q$-gons, which is glued, by the boundary circle, with a
disc of (the half of original) $p$-gons.

\break
\begin{center}
\d{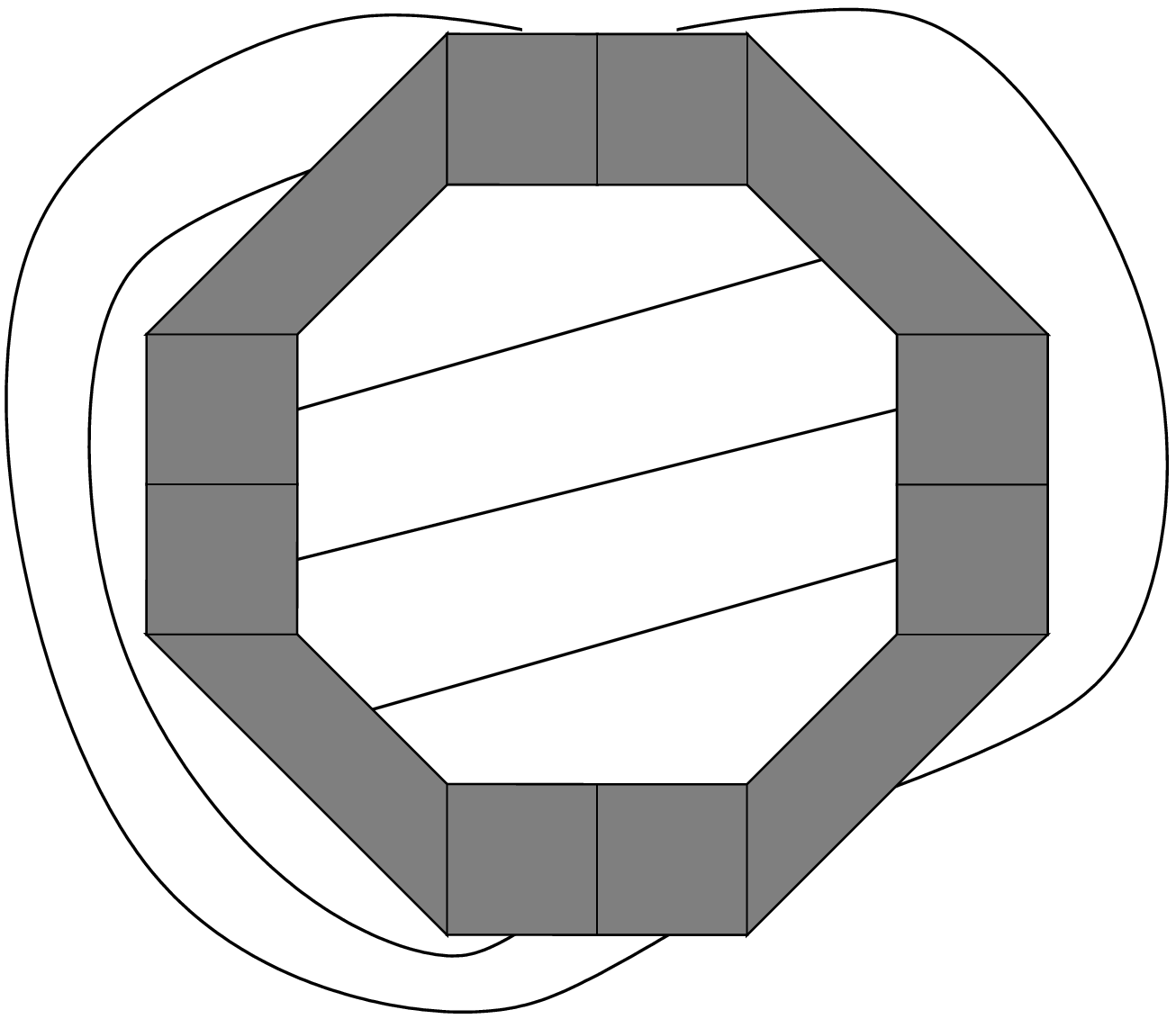}\hspace{4cm}\d{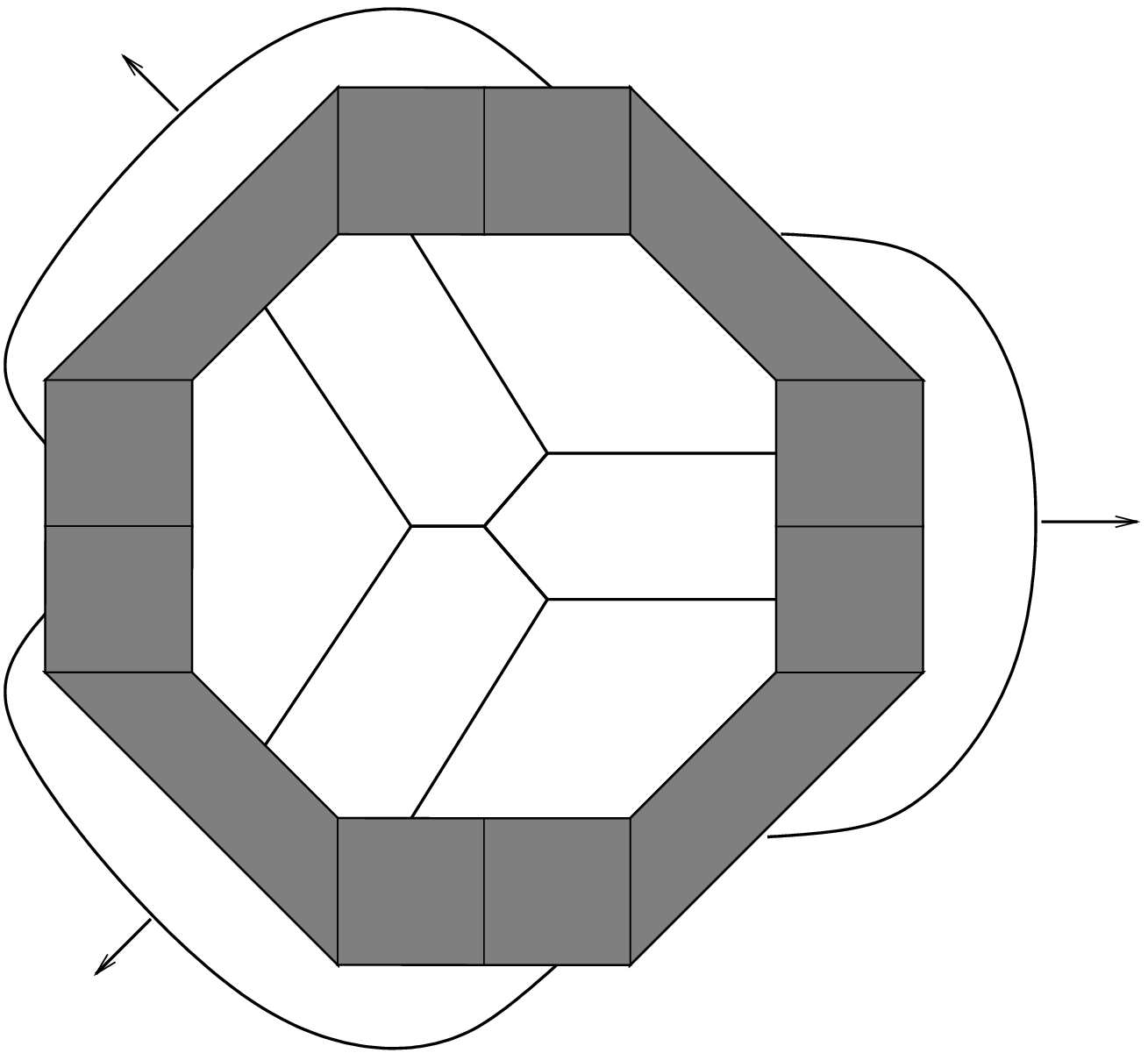}\\[3.5mm]

\d{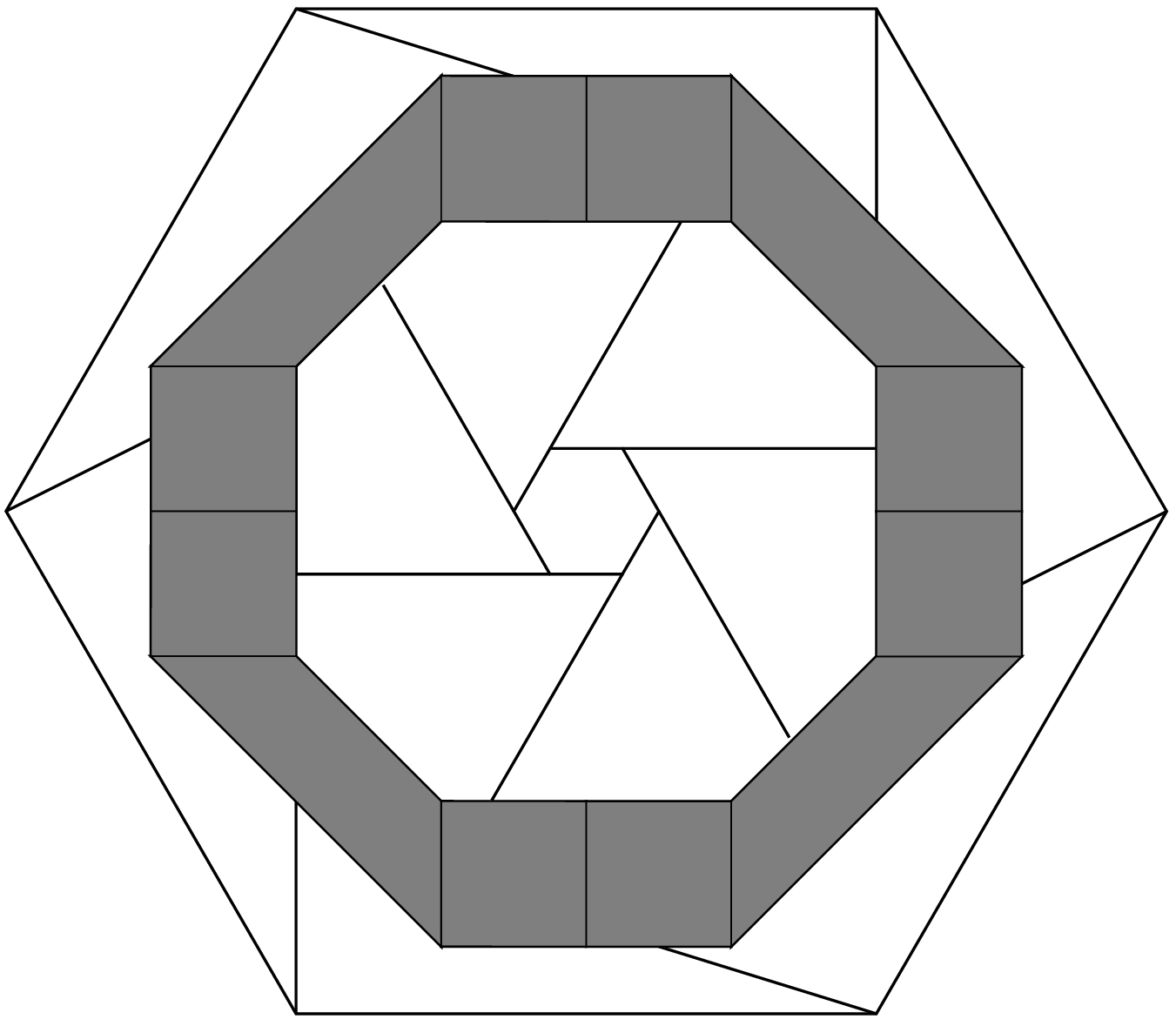}\hspace{4cm}\d{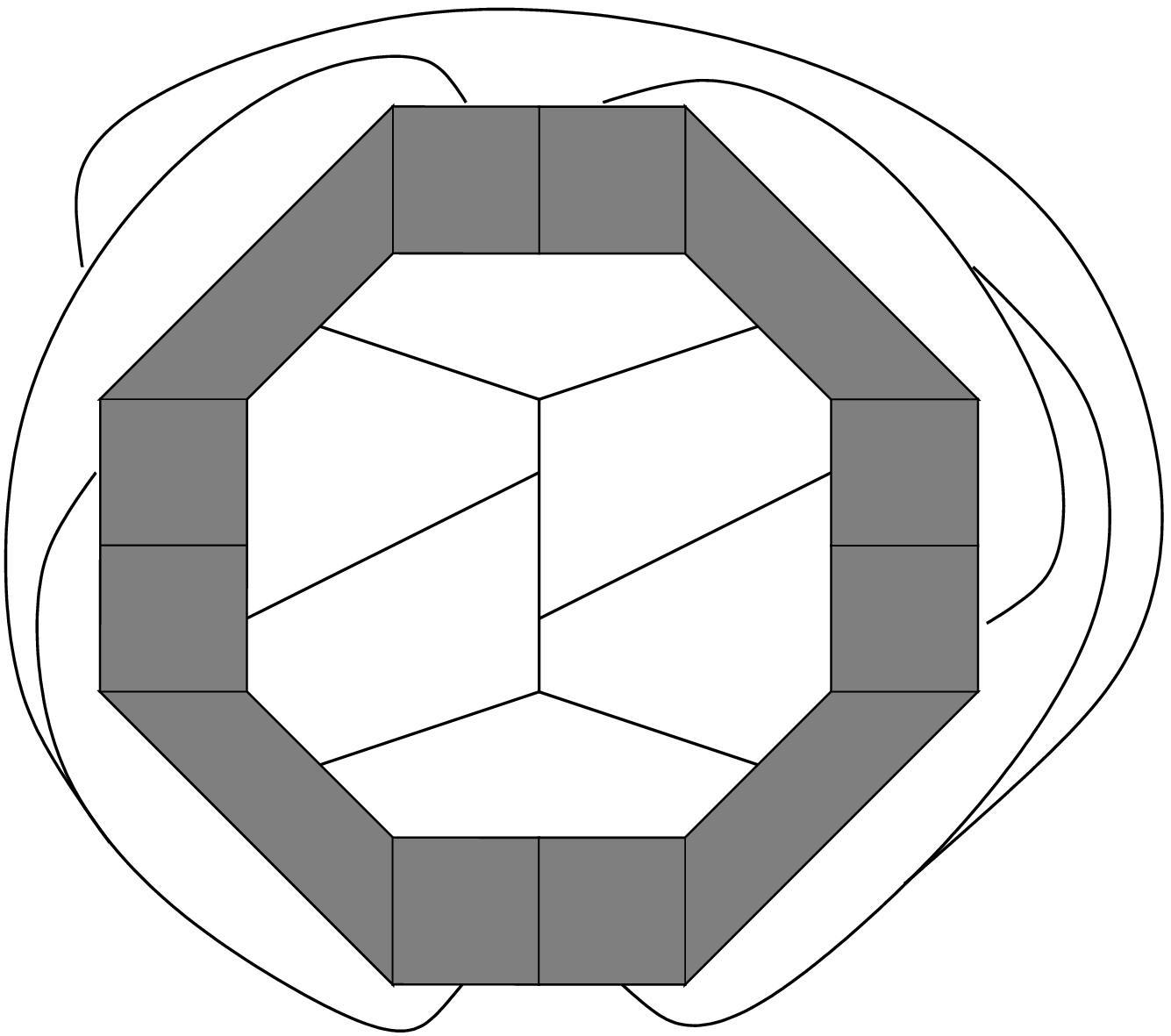}\\[3.5mm]

\d{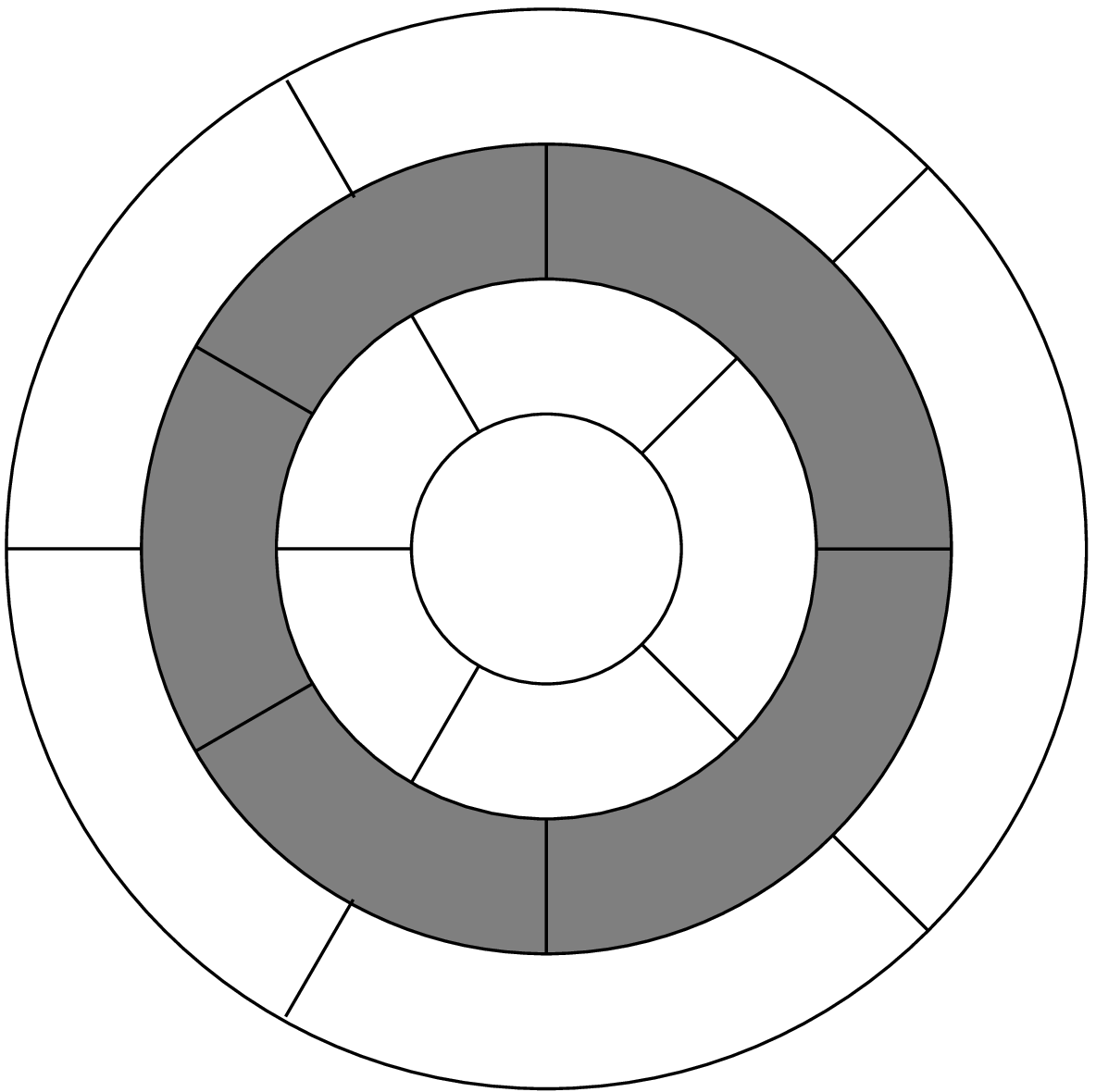}\hspace{4cm}\d{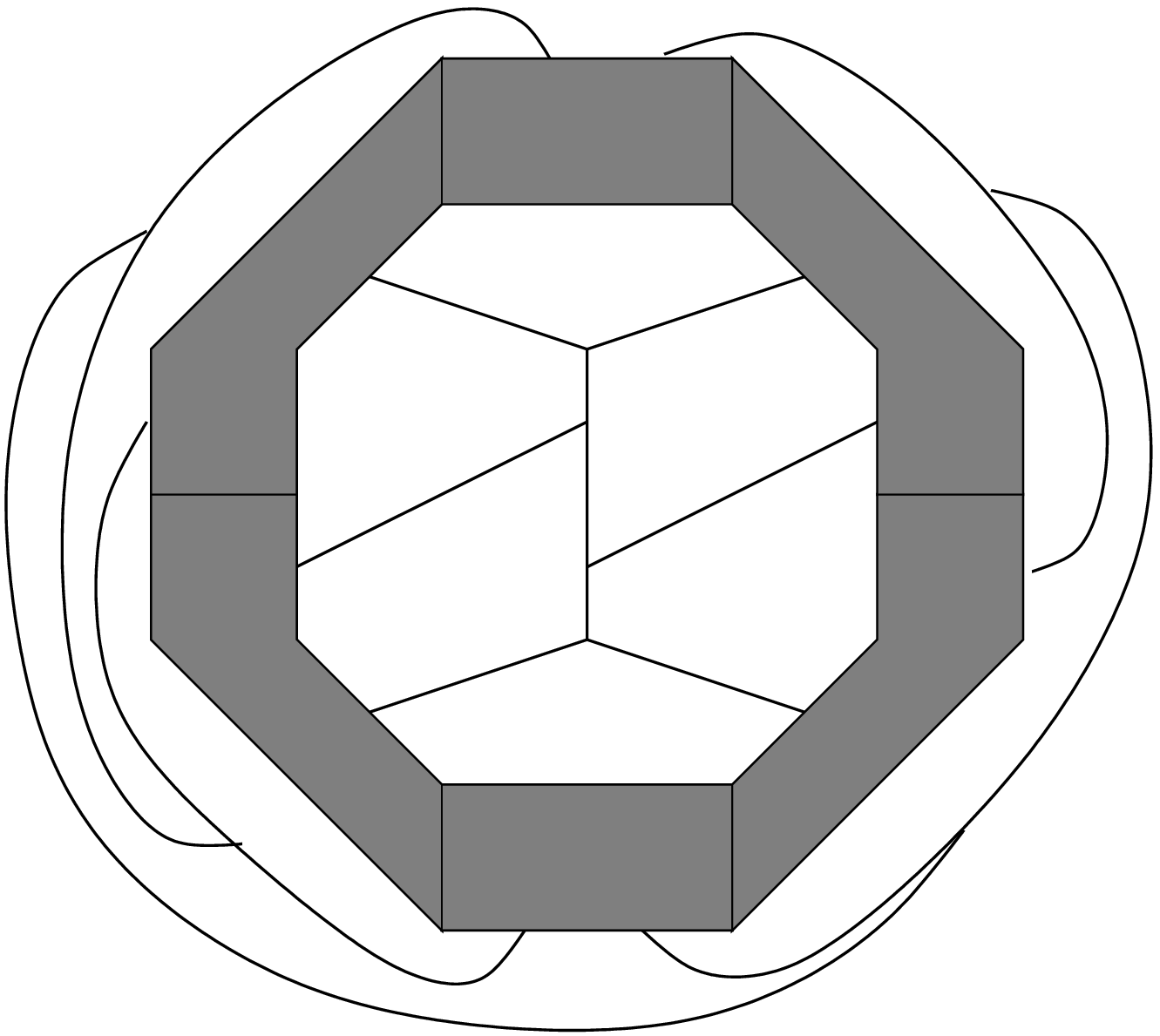}\\[3.5mm]

\d{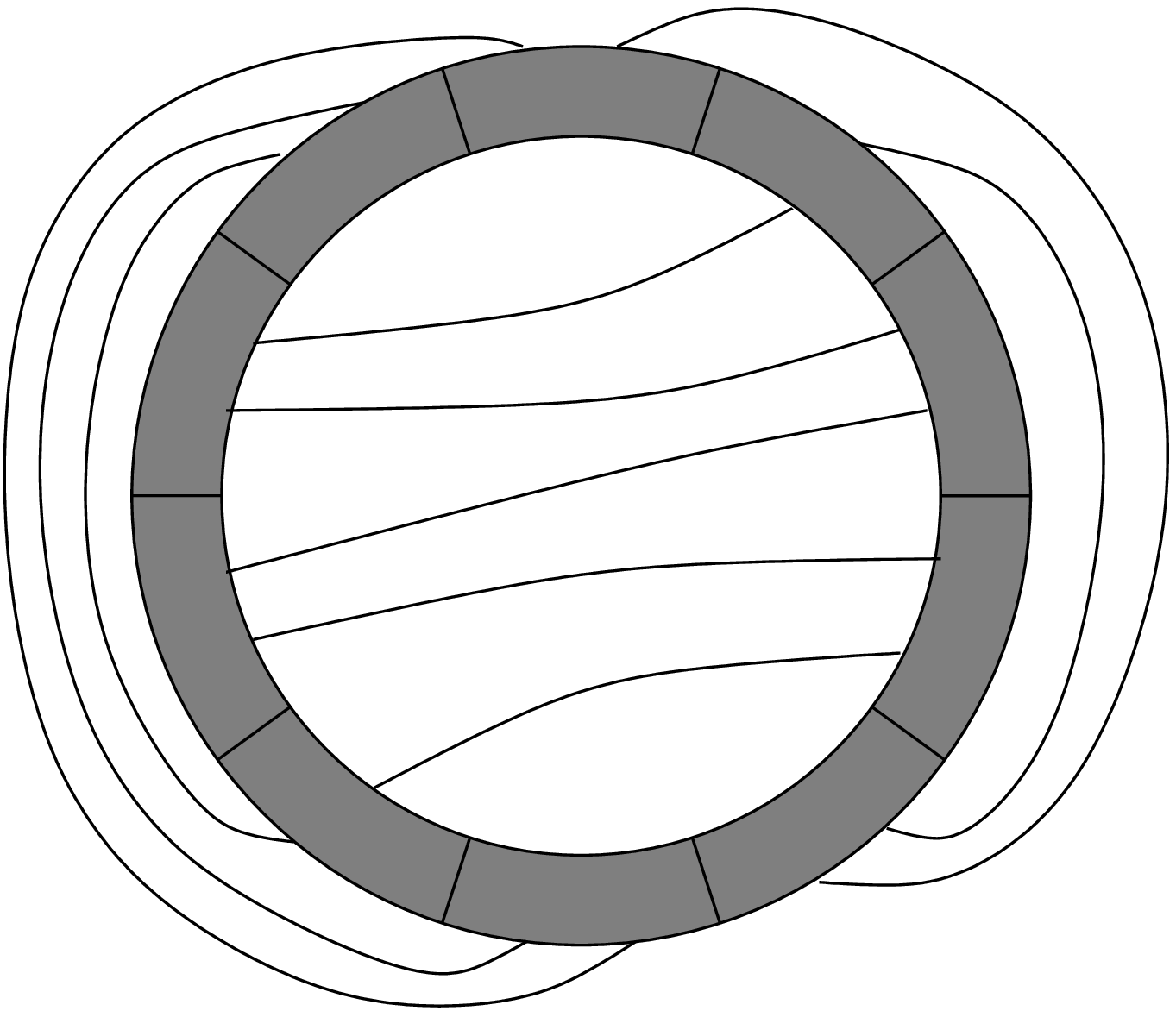}\hspace{4cm}\d{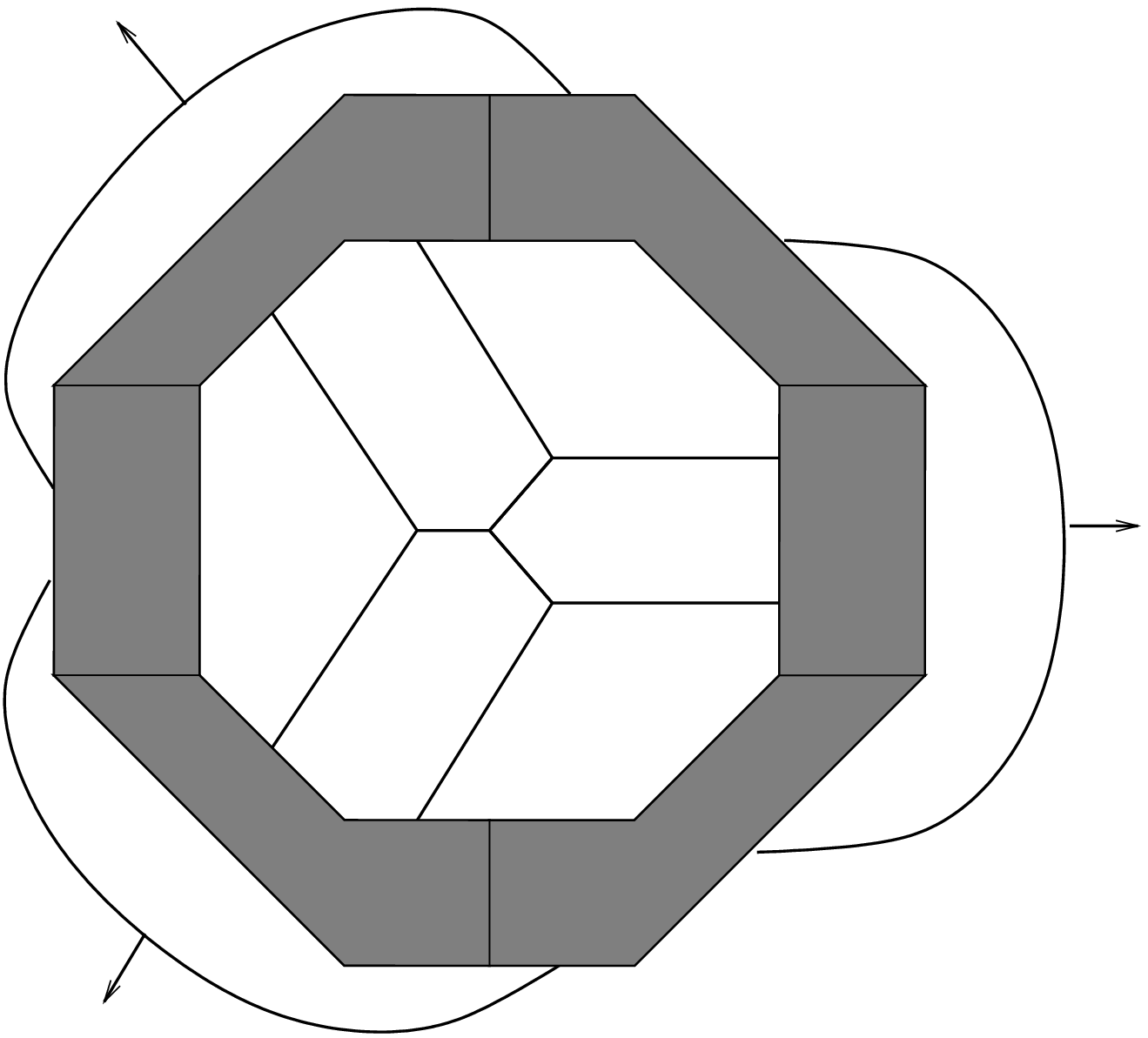}\\[3.5mm]

\d{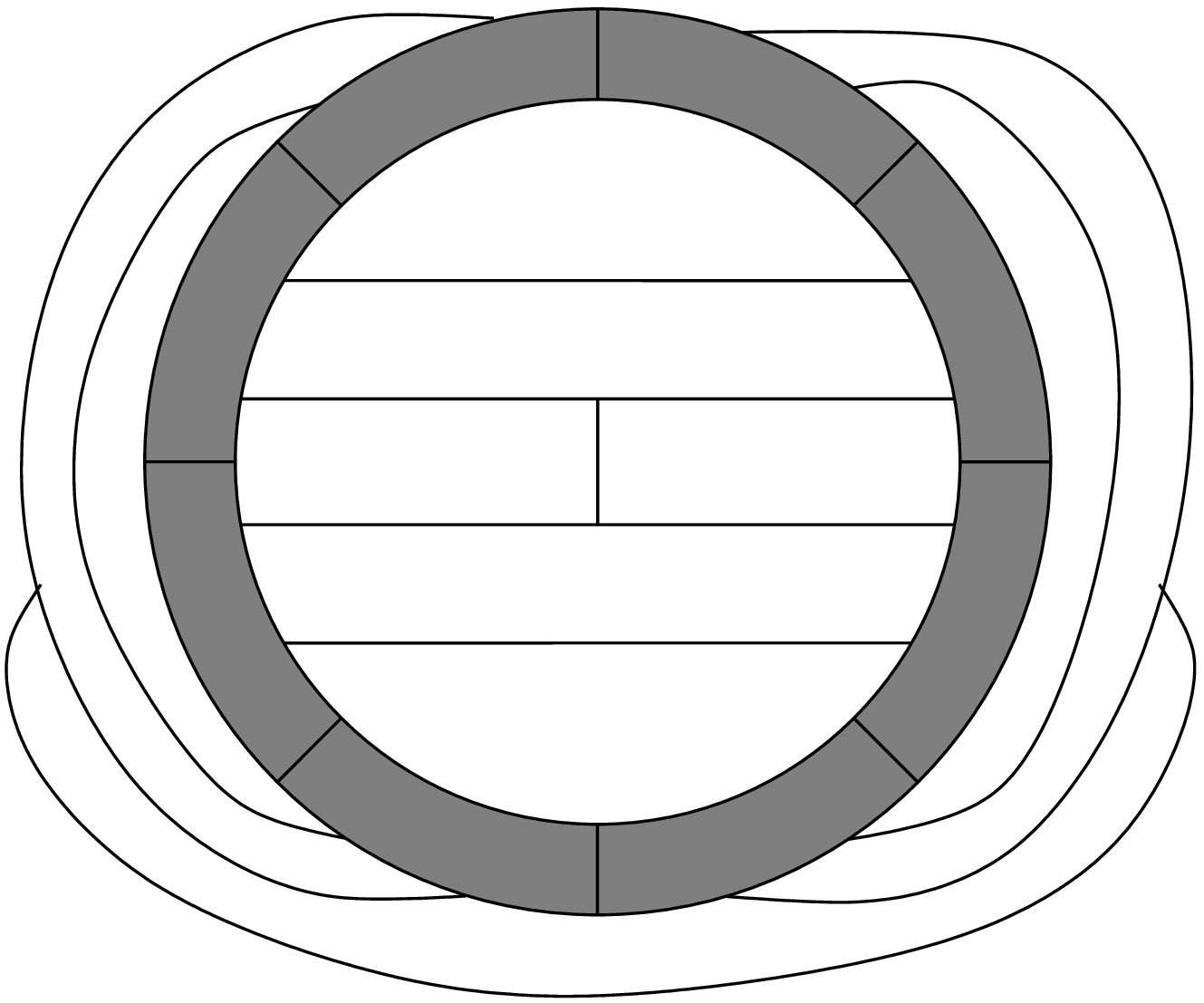} \\[3mm]

% \vspace{5mm}
Figure 1: {\em All fullerene maps $M_n(p,q)$ (i.e. with $ \{ p,q \} = \{ 5,6 \}
$): four maps $M_{12}(6,5)$ and five maps $M_n(5,6)$ with $n=5,6,10,6,8$.}
\end{center}

% \vspace{5mm}

\begin{center}
\dd{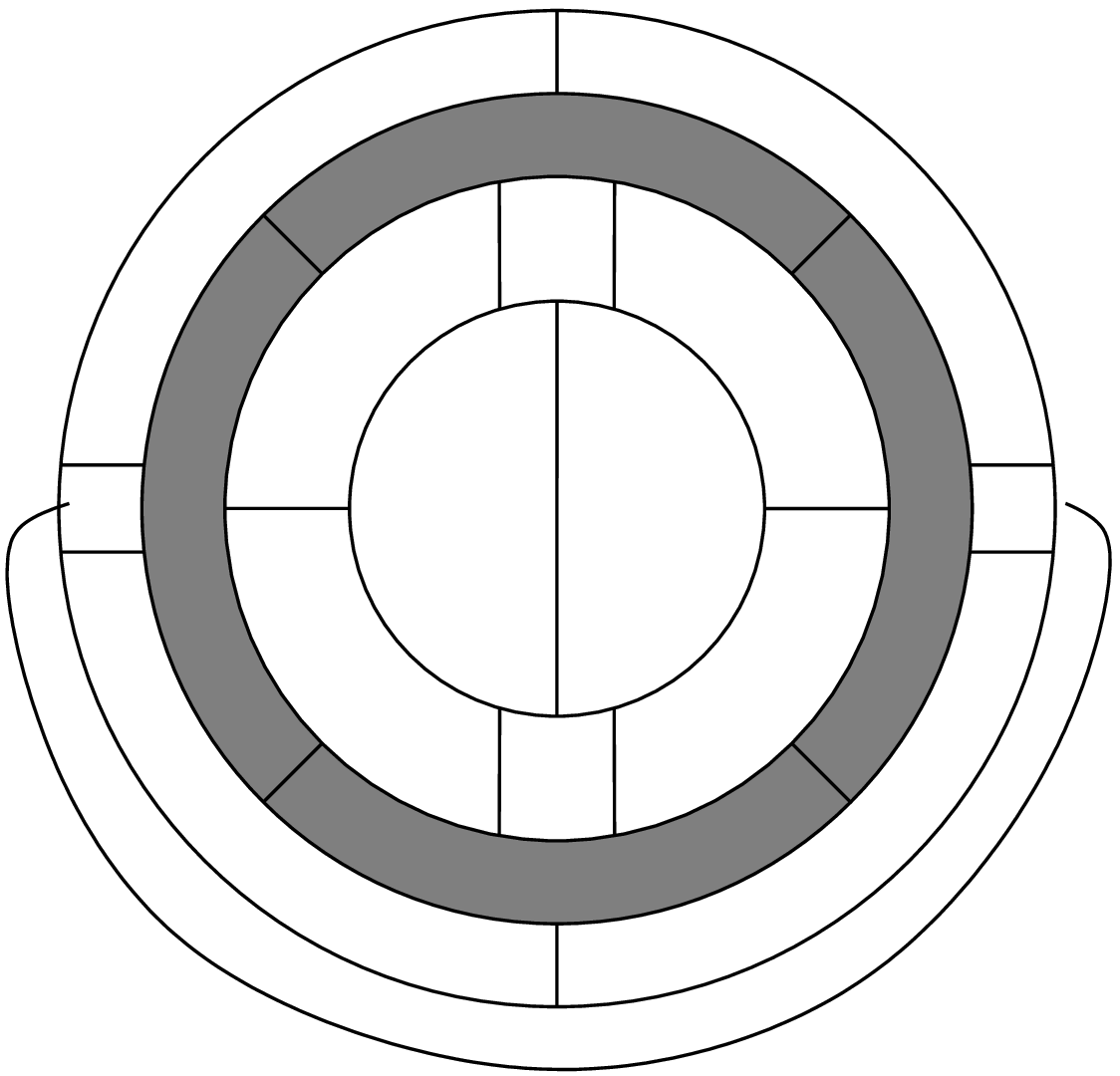}\hspace{3.7cm}\da{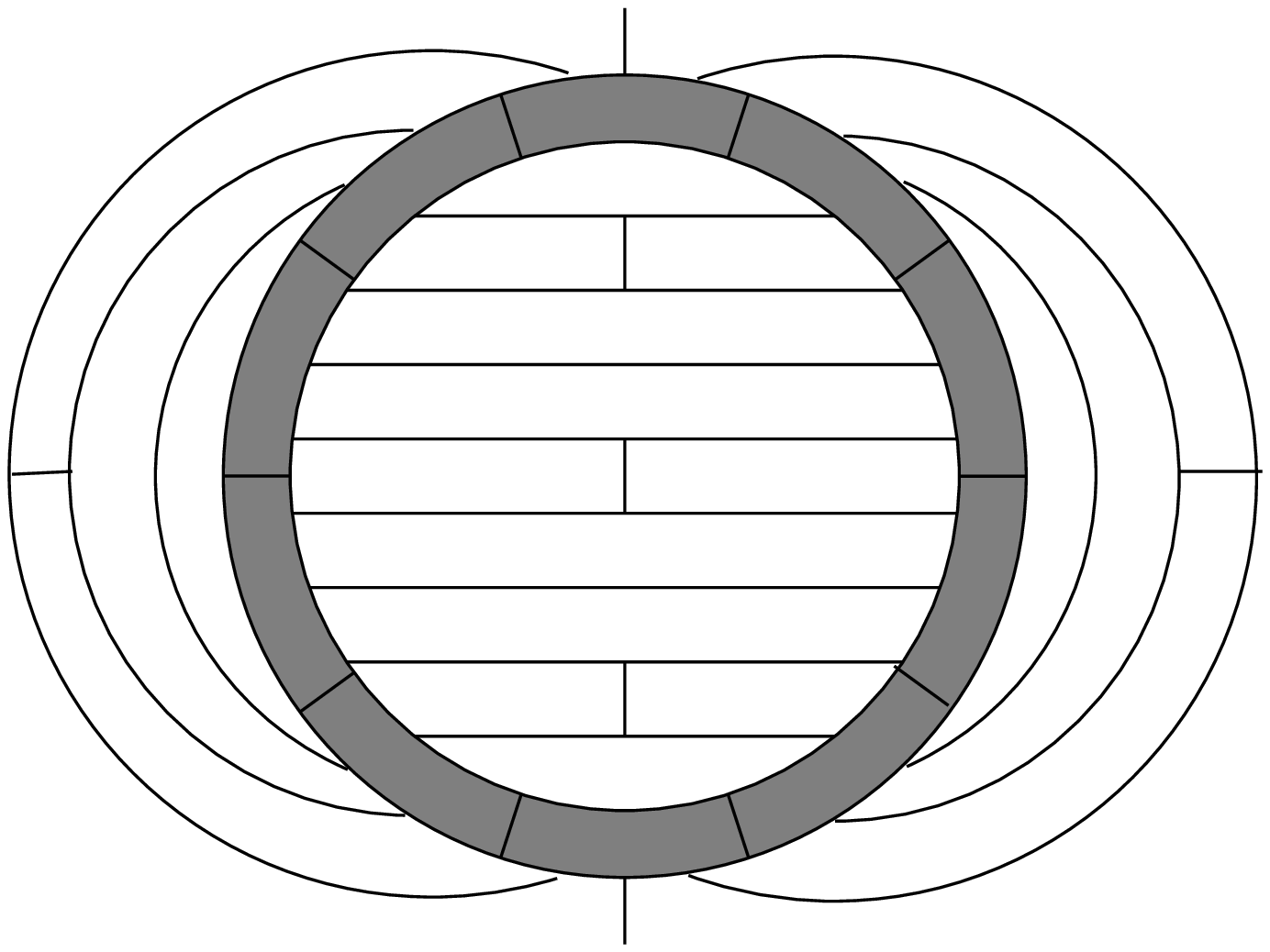}\\[3.3mm]

\d{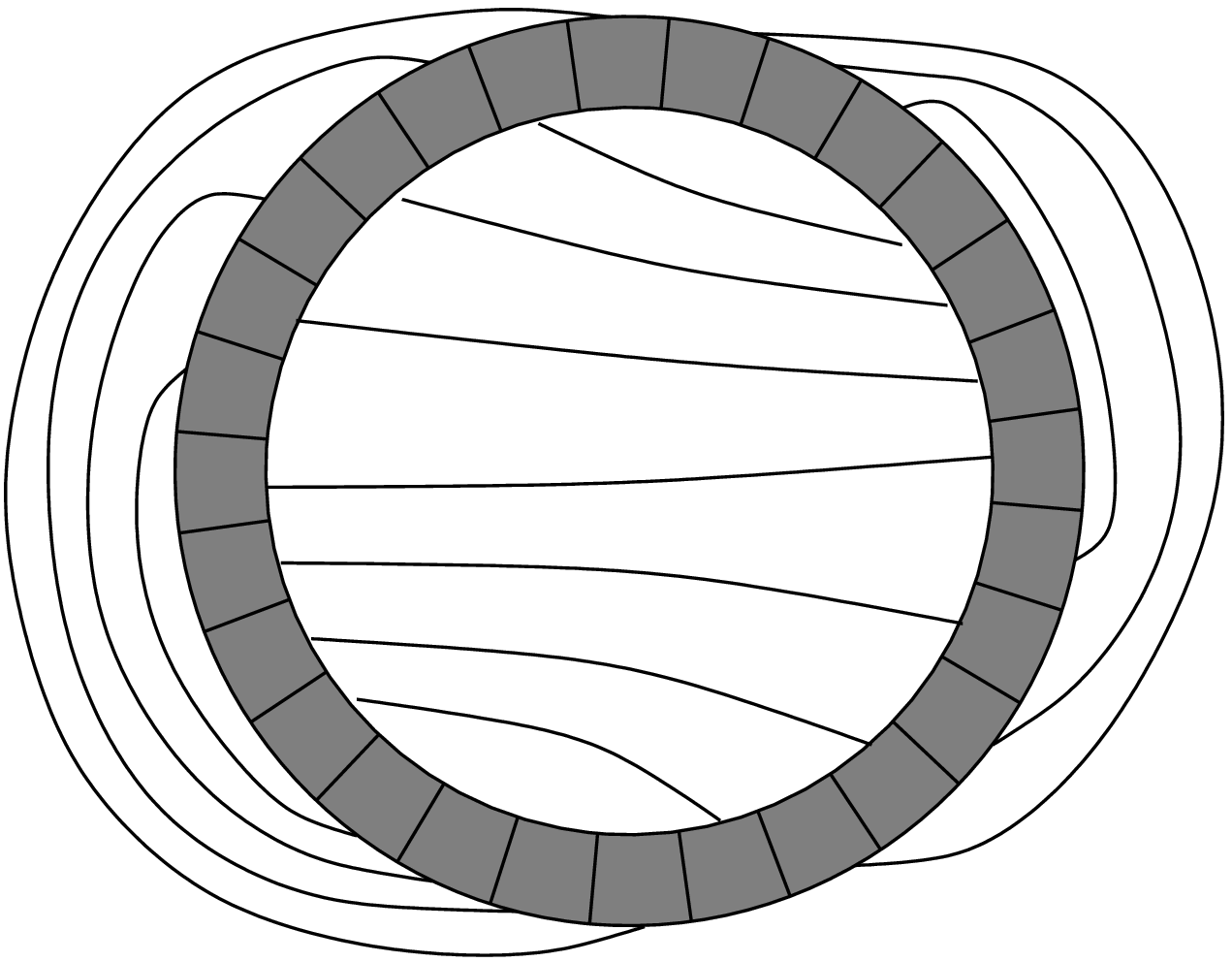}\hspace{4cm}\dd{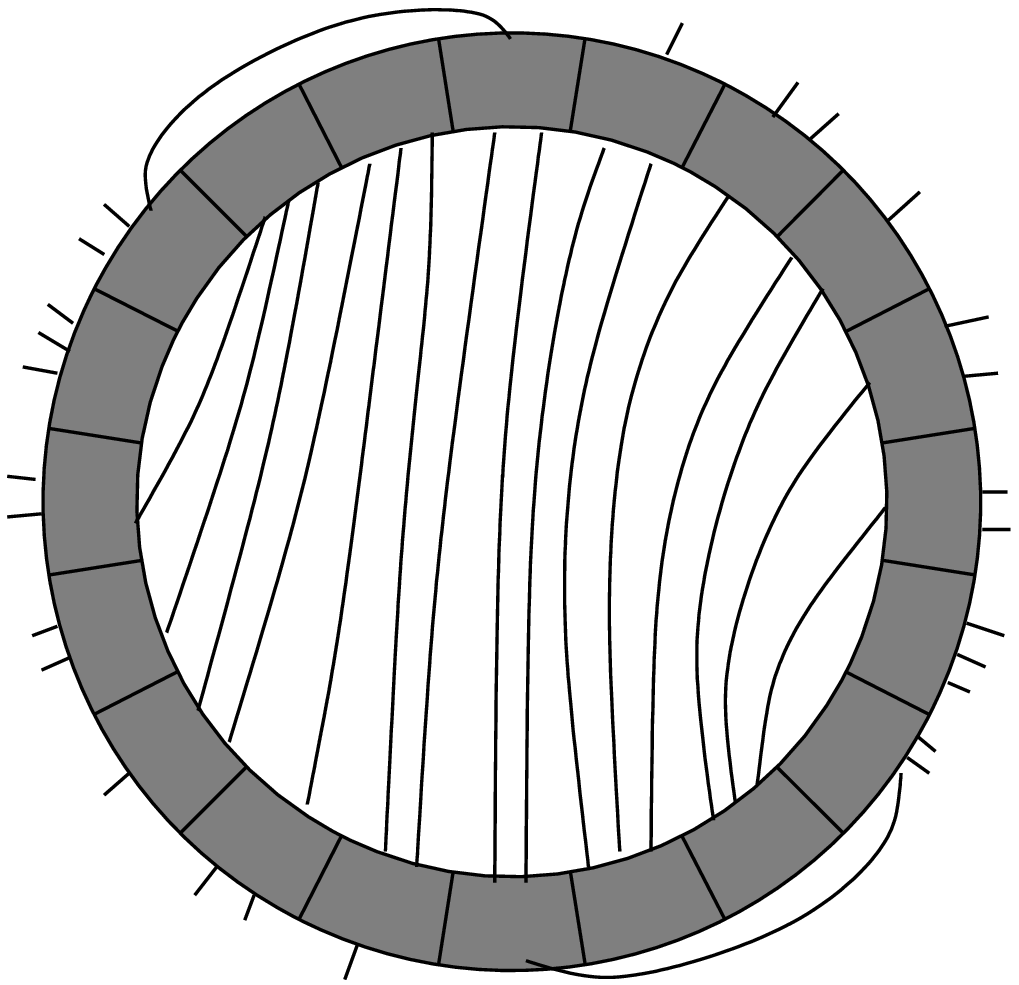}\\[3.3mm]

\dd{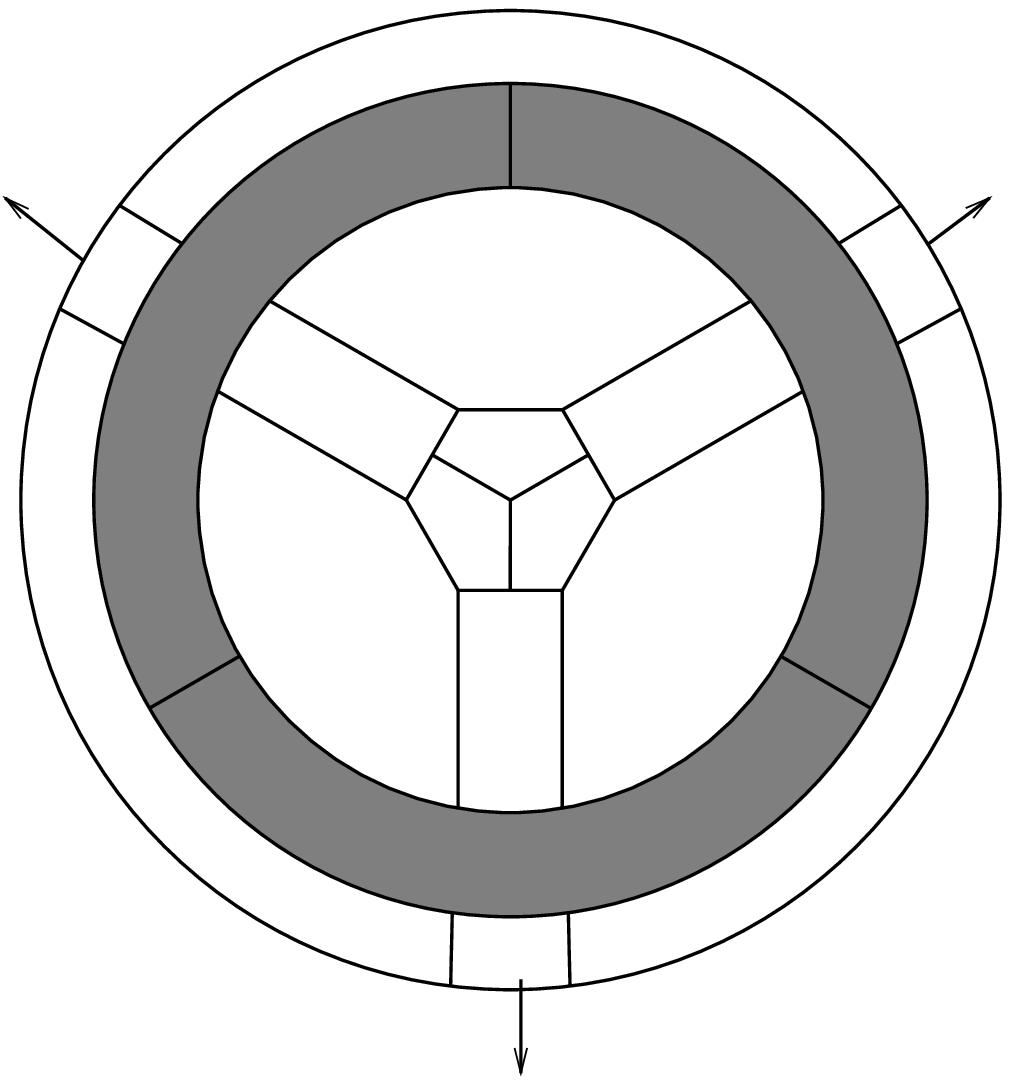}\hspace{4cm}\dd{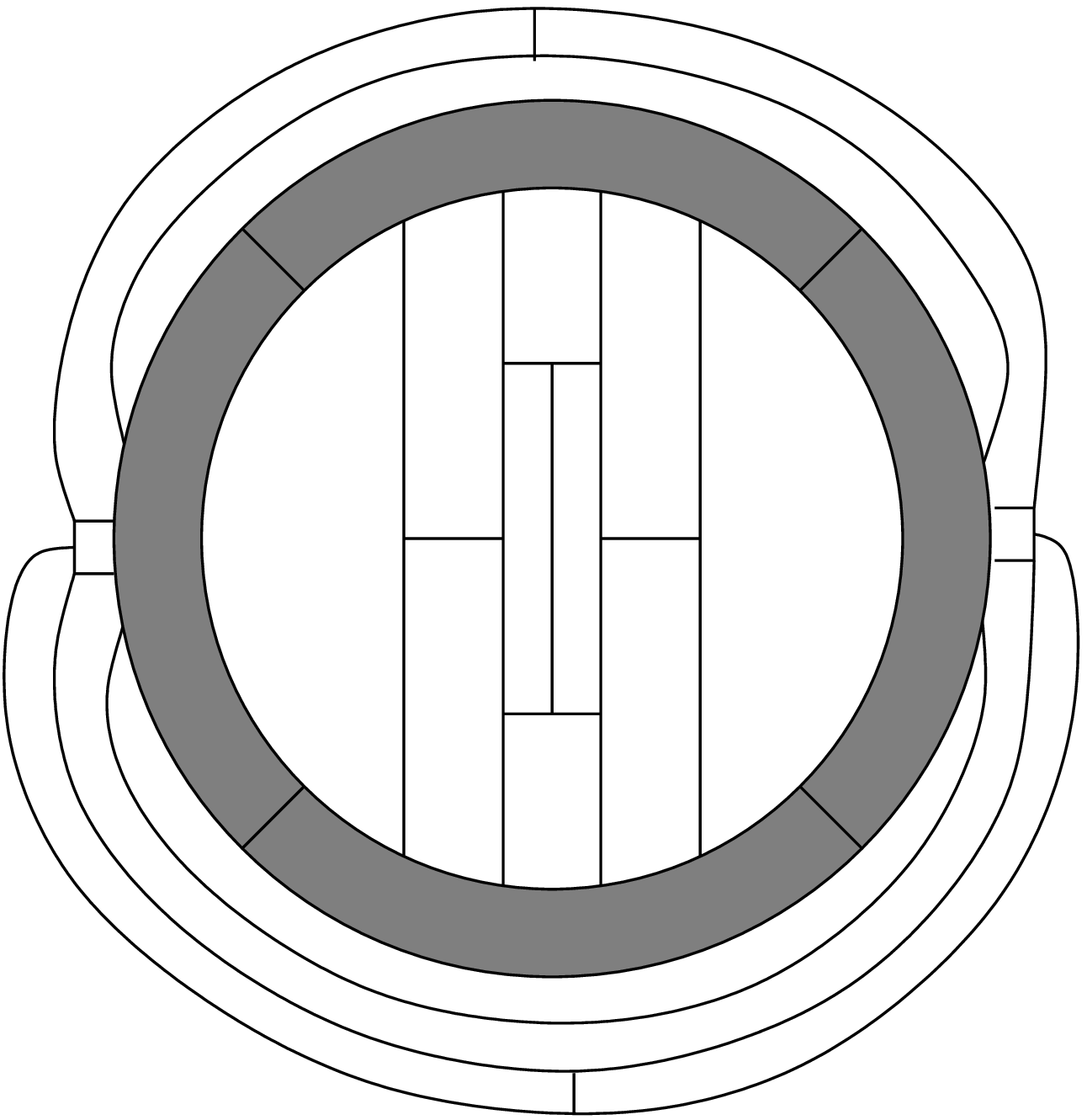}\\[3.3mm]

\dd{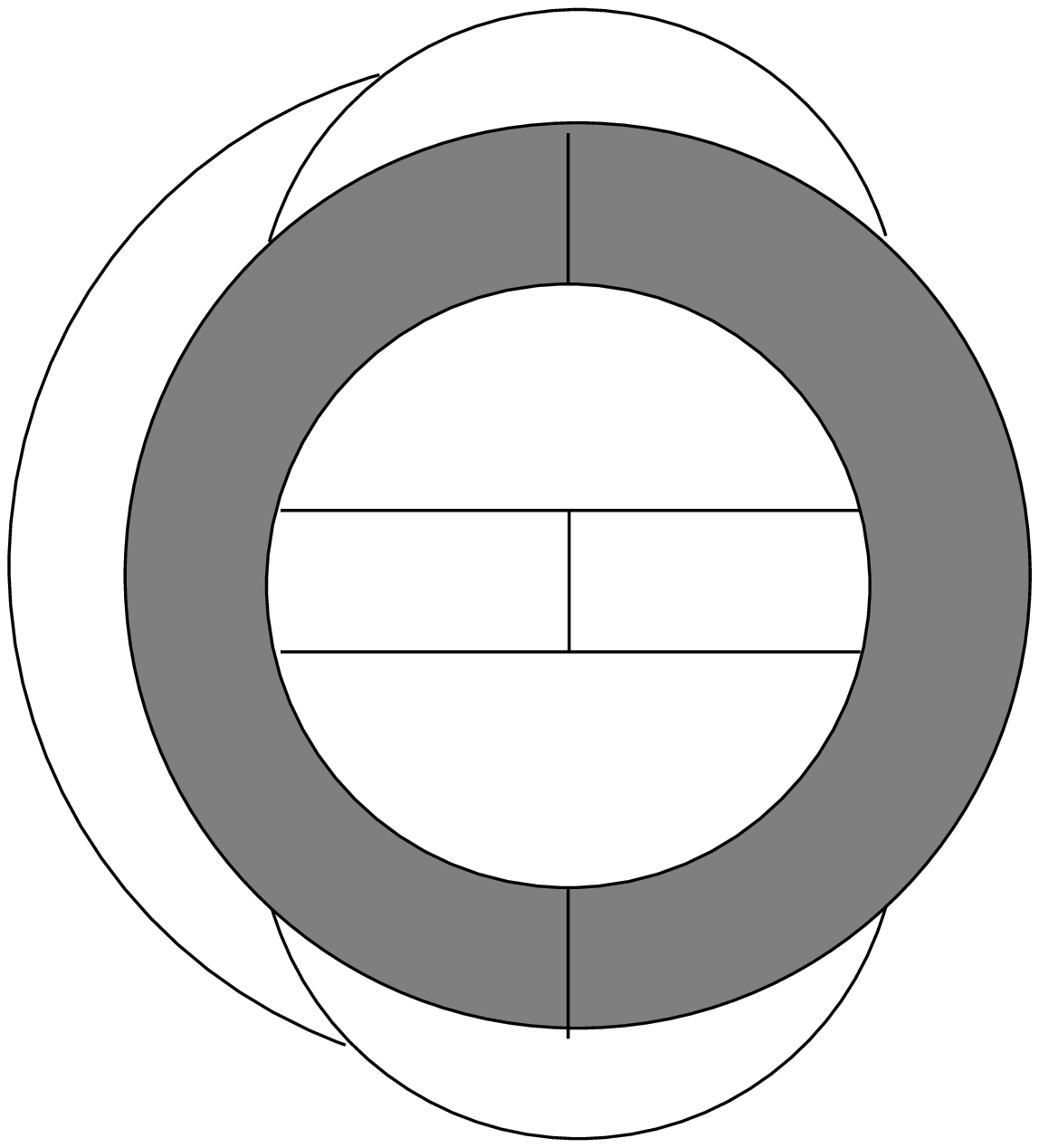}\hspace{4cm}\d{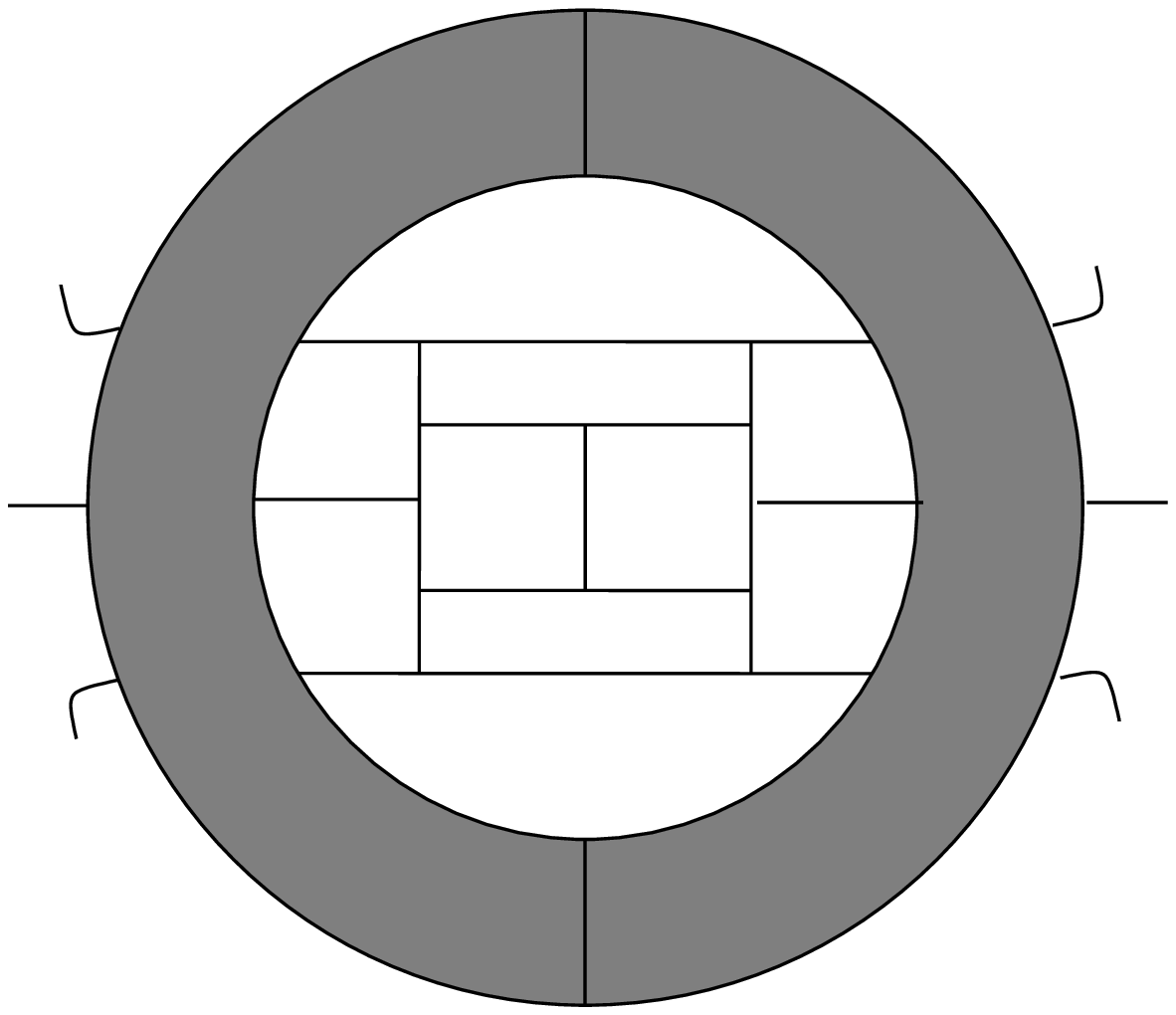}\\[3.3mm]

\ddd{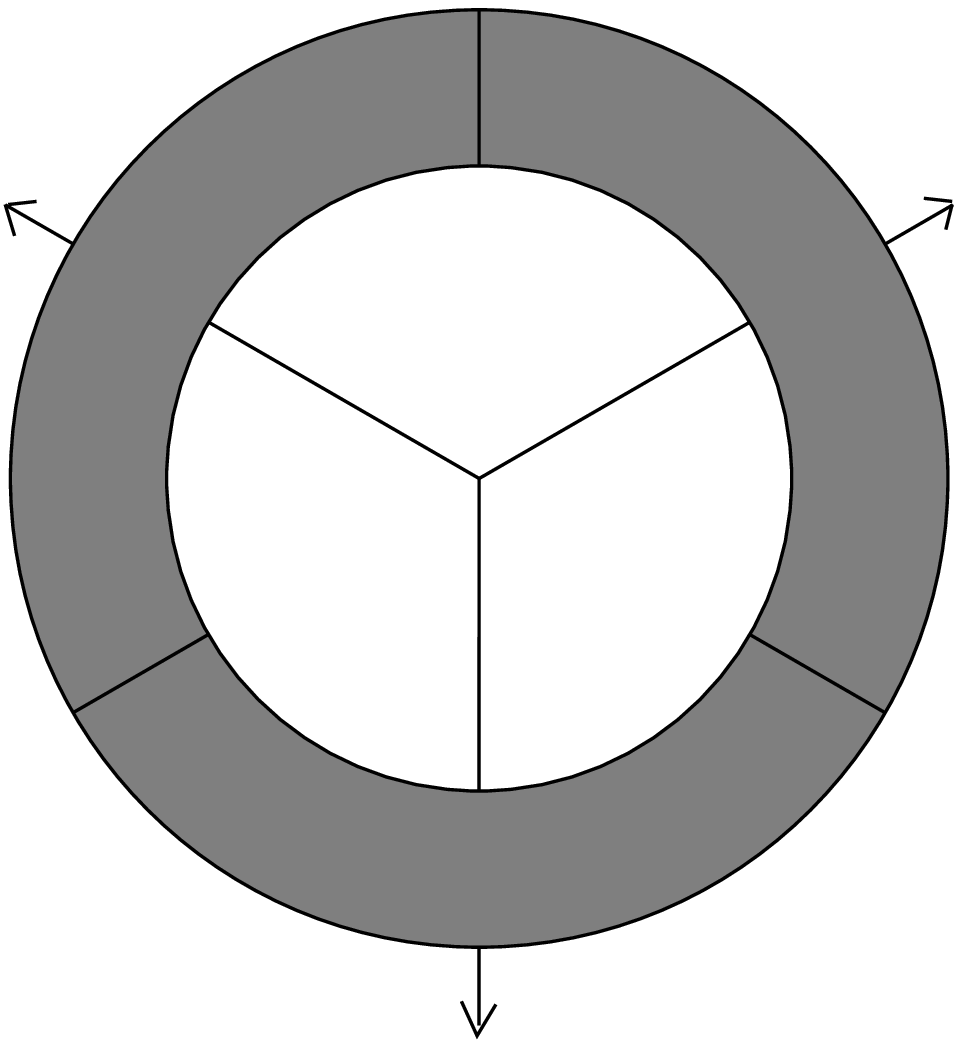}\hspace{4cm}\ddd{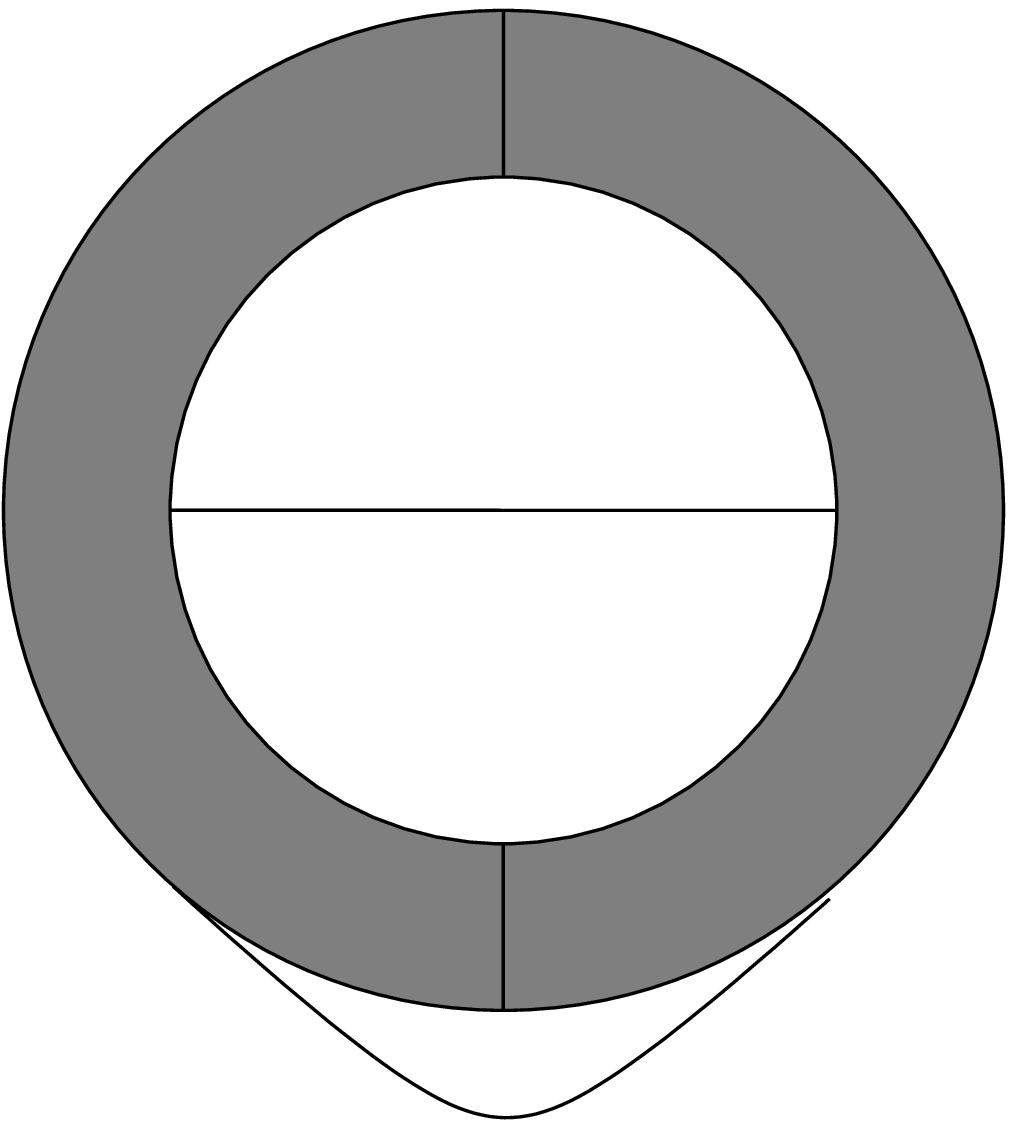}\\[3.3mm]

Figure 2: {\em Four azulenoids ($M_4(5,7)$, $M_{10}(5,7)$, $M_{28}(7,5)$,
$M_{20}(5,7)$), two $(5,8)$-maps ($M_3(5,8)$, $M_4(5,8)$) and 
$M_2(4,8)$, $M_2(5,10)$, $M_3(4,6)$, $M_2(3,6)$.}
\end{center}

\setlength{\unitlength}{1cm}
\begin{minipage}[t]{3.5cm}
\begin{picture}(3.5,3.5)
\leavevmode
\epsfxsize=3.5cm
\epsffile{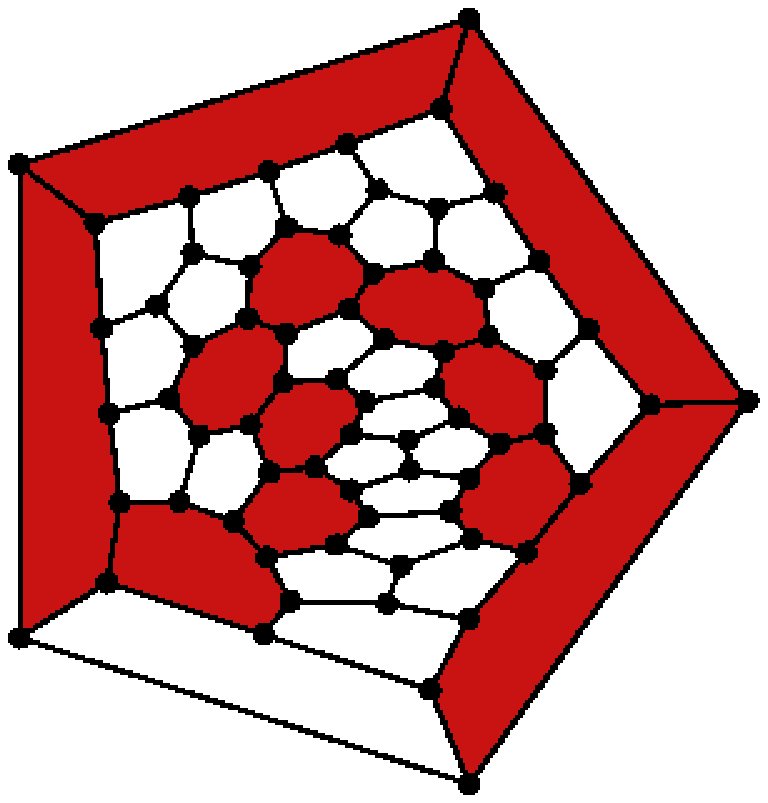}
\end{picture}\par
\begin{center}
%{{\bf Nr.1} \quad \\}
\end{center}
\end{minipage}
\setlength{\unitlength}{1cm}
\begin{minipage}[t]{3.5cm}
\begin{picture}(3.5,3.5)
\leavevmode
\epsfxsize=3.5cm
\epsffile{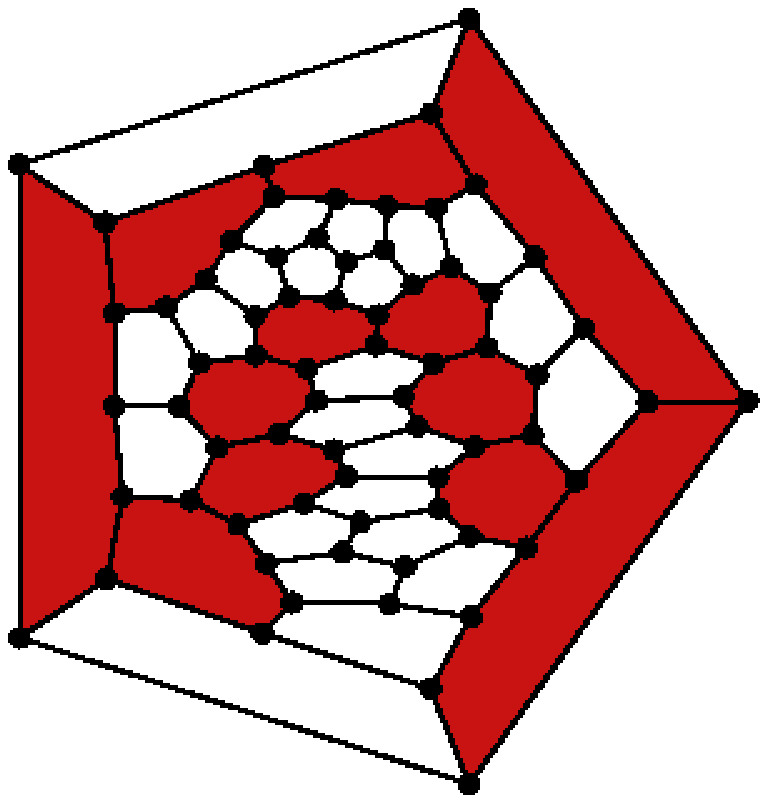}
\end{picture}\par
\begin{center}
%{{\bf Nr.2} \quad \\}
\end{center}
\end{minipage}
\setlength{\unitlength}{1cm}
\begin{minipage}[t]{3.5cm}
\begin{picture}(3.5,3.5)
\leavevmode
\epsfxsize=3.5cm
\epsffile{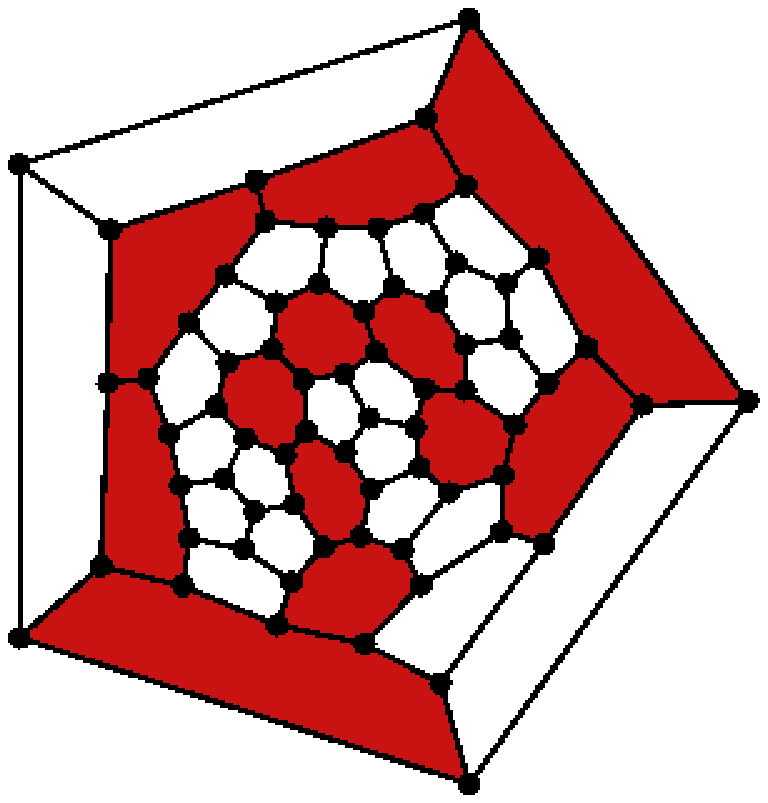}
\end{picture}\par
\begin{center}
%{{\bf Nr.3} \quad \\}
\end{center}
\end{minipage}

\setlength{\unitlength}{1cm}
\begin{minipage}[t]{3.5cm}
\begin{picture}(3.5,3.5)
\leavevmode
\epsfxsize=3.5cm
\epsffile{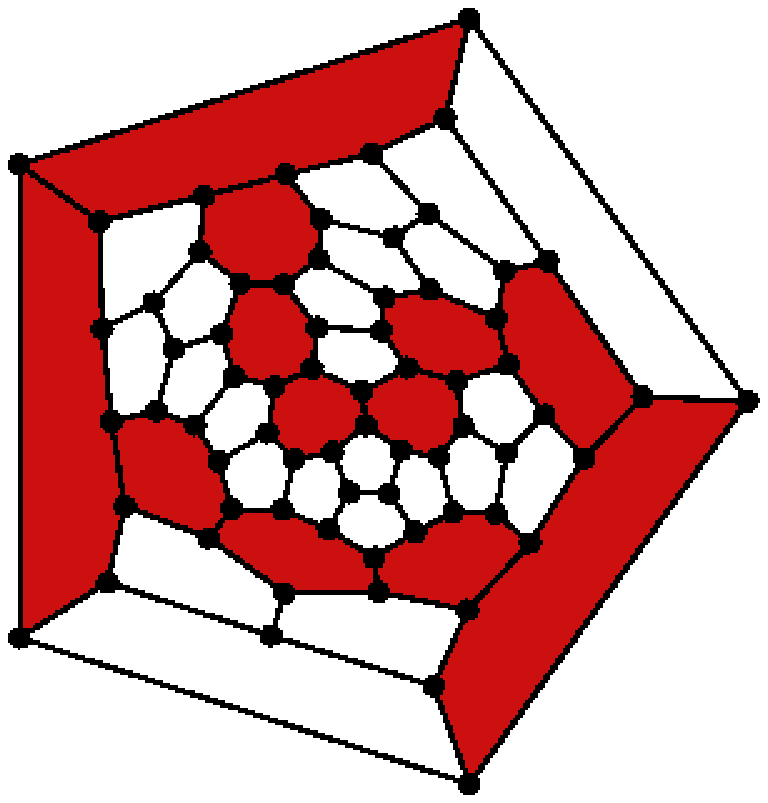}
\end{picture}\par
\begin{center}
%{{\bf Nr.4} \quad \\}
\end{center}
\end{minipage}
\setlength{\unitlength}{1cm}
\begin{minipage}[t]{3.5cm}
\begin{picture}(3.5,3.5)
\leavevmode
\epsfxsize=3.5cm
\epsffile{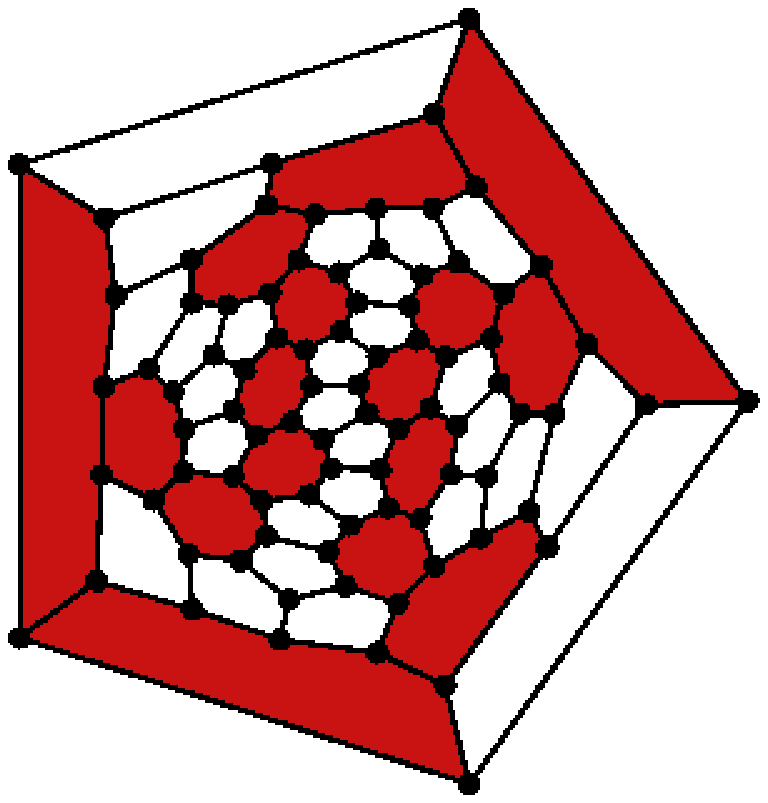}
\end{picture}\par
\begin{center}
%{{\bf Nr.5} \quad \\}
\end{center}
\end{minipage}
\setlength{\unitlength}{1cm}
\begin{minipage}[t]{3.5cm}
\hfil\begin{picture}(3.5,3.5)
\leavevmode
\epsfxsize=3.5cm
\epsffile{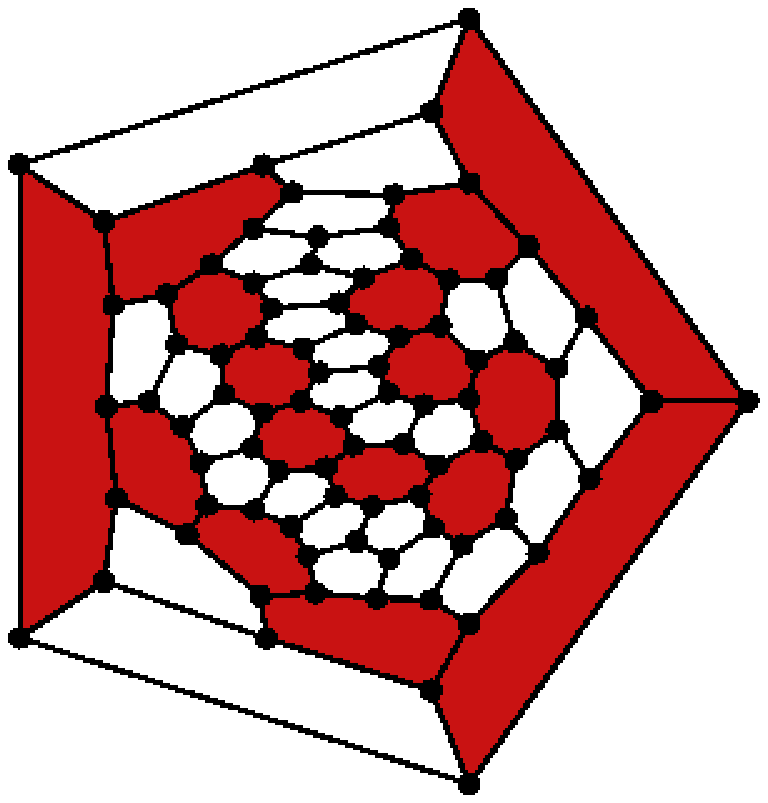}
\end{picture}\hfil\par
\begin{center}
%{{\bf Nr.6} \quad \\}
\end{center}
\end{minipage}

        \begin{center}
	Figure 3: azulenoids $M_{n}(5,7)$ from \cite{Ha}: 4 with
$n=12$ and 2 with $n=16$.
	\end{center}

\end{document}